\documentclass[11pt]{elsarticle}

\makeatletter
\def\ps@pprintTitle{%
 \let\@oddhead\@empty
 \let\@evenhead\@empty
 \def\@oddfoot{}%
 \let\@evenfoot\@oddfoot}

\usepackage[left=2.5cm,right=2.5cm,top=2cm,bottom=3cm]{geometry}
\usepackage{algorithm}
\usepackage{algorithmic}
\usepackage{subfigure} 
\usepackage{amsfonts,amssymb,amsmath,amsthm,bm}
\numberwithin{equation}{section}
\usepackage{enumerate}
\usepackage{multirow}
\usepackage{xcolor}
\usepackage{graphicx,float,epstopdf}
\usepackage{array}
\usepackage{booktabs}
\usepackage{hyperref}
\usepackage{subeqnarray}
\usepackage{cases}
\newcommand\highlightReference[1]{%
  \expandafter\newcommand\csname highlightReference-#1\endcsname{}%
}
\let\oldbibitem\bibitem
\def\bibitem#1 #2\par{%
  \expandafter\ifx\csname highlightReference-#1\endcsname\relax
    \oldbibitem{#1}#2\par
  \else
    \oldbibitem{#1}\highlight{#2}\par
  \fi
}
\usepackage{color}
\newcommand\highlight[1]{\textcolor{blue}{#1}}

\newtheorem{theorem}{Theorem}[section]

\newtheorem{property}[theorem]{Property}

\newtheorem{remark}[theorem]{Remark}

\theoremstyle{definition}
\newtheorem{example}{Example}[section]

\begin{document}

\title{{\bfseries Second-order flows for computing the ground states of rotating Bose-Einstein condensates}}

\author[hunnu]{Haifan Chen}\ead{haifanchen@smail.hunnu.edu.cn}

\author[csu]{Guozhi Dong\corref{cor}}\ead{guozhi.dong@csu.edu.cn}
\cortext[cor]{Corresponding author.}

\author[scnu]{Wei Liu}\ead{wliu@m.scnu.edu.cn}

\author[hunnu]{Ziqing Xie}\ead{ziqingxie@hunnu.edu.cn}

\address[hunnu]{Key Laboratory of Computing and Stochastic Mathematics (Ministry of Education), School of Mathematics and Statistics, Hunan Normal University, Changsha 410081, China}
\address[csu]{School of Mathematics and Statistics, HNP-LAMA,  Central South University, Changsha 410083, China}
\address[scnu]{South China Research Center for Applied Mathematics and Interdisciplinary Studies, South China Normal University, Guangzhou 510631, China}

\date{\today}


\begin{frontmatter}

\begin{abstract}
Second-order flows in this paper refer to some artificial evolutionary differential equations involving second-order time derivatives distinguished from gradient flows which are considered to be first-order flows. This is a popular topic due to the recent advances of inertial dynamics with damping in convex optimization. Mathematically, the ground state of a rotating Bose-Einstein condensate (BEC) can be modeled as a minimizer of the Gross-Pitaevskii energy functional with angular momentum rotational term under the normalization constraint. We introduce two types of second-order flows as energy minimization strategies for this constrained non-convex optimization problem, in order to approach the ground state. The proposed artificial dynamics are novel second-order nonlinear hyperbolic partial differential equations with dissipation. Several numerical discretization schemes are discussed, including explicit and semi-implicit methods for temporal discretization, combined with a Fourier pseudospectral method for spatial discretization. These provide us a series of efficient and robust algorithms for computing the ground states of rotating BECs. Particularly, the newly developed algorithms turn out to be superior to the state-of-the-art numerical methods based on the gradient flow. {In comparison with the gradient flow type approaches: When explicit temporal discretization strategies are adopted, the proposed methods allow for larger stable time step sizes; While for semi-implicit discretization, using the same step size, a much smaller number of iterations are needed for the proposed methods to reach the stopping criterion, and every time step encounters almost the same computational complexity.} 
Rich and detailed numerical examples are documented for verification and comparison.
\end{abstract}

\begin{keyword}
rotating Bose-Einstein condensate, 
Gross-Pitaevskii energy functional,
ground state, 
inertial dynamics, 
second-order dissipative hyperbolic PDEs, 
constrained non-convex minimization
\end{keyword}

\end{frontmatter}


\section{Introduction}

\subsection{Context and background}

The renowned Bose-Einstein condensate (BEC) is a state of matter that was predicted theoretically by S. N. Bose and A. Einstein in 1924–1925 and first realized experimentally in 1995 in dilute bosonic atomic gases at extremely low temperatures \cite{AEMWC1995Science, BSTH1995PRL, DMADDKK1995PRL}. Up to now, BEC has been created in a variety of other bosonic systems (e.g., molecules, quasiparticles and photons) and in different environments (e.g., in an Earth-orbiting research lab), and has become an underlying tool/platform for experiments in many related fields such as superfluids, supersolids, quantum information processing and precision atomic devices \cite{BECspace,FS2001JPCM,klaers2010bose,supersolidity}. As one of the most active directions of BECs, the study of quantized vortices has attracted great interests from both the physics and mathematics communities, and in particular it provides valuable perspectives for exploring mysterious superfluids \cite{BWM2005CMS,FS2001JPCM}. Experimentally, quantized vortices could usually be generated from the ground state of rotating BECs \cite{ARVK2001Science,ECHC2002PRL,Cornell2003giant,MAHHWC1999PRL,MCWD2000PRL}. Meanwhile, with the growth of computing power, numerical simulation to efficiently find the ground state can facilitate a better understanding of theories and vortex properties of rotating BECs as well as predicting and guiding experiments.

Within the mean-field theory, a rotating BEC at zero or very low temperature can be modeled by the Gross-Pitaevskii (GP) equation with an angular momentum rotational term \cite{BC2013KRM,BDZ2006SIAP,DGPS1999RMP,Fetter2009RMP,Seiringer2002CMP}. In this setting, after proper nondimensionalization and dimension reduction, the ground state is defined as a complex-valued macroscopic wave function (or order parameter) $\phi_g:\mathbb{R}^d\to\mathbb{C}$ ($d=2,3$) that minimizes the GP energy functional \cite{BC2013KRM,BWM2005CMS,DGPS1999RMP,LSY2000PRA,Seiringer2002CMP}
\begin{align}\label{eq:GPErot-energy}
E(\phi) = \int_{\mathbb{R}^d} \left( \frac12|\nabla\phi(\mathbf{x})|^2+V(\mathbf{x})|\phi(\mathbf{x})|^2+\frac{\beta}{2}|\phi(\mathbf{x})|^4 -\Omega\overline{\phi}(\mathbf{x})L_z\phi(\mathbf{x})\right) \mathrm{d}\mathbf{x}
    \end{align}
under the normalization constraint
\begin{align}\label{eq:ELeq-cons}
\|\phi\|^2:=\int_{\mathbb{R}^d}|\phi(\mathbf{x})|^2\mathrm{d}\mathbf{x}=1.
\end{align}
Here, $\mathbf{x}=(x,y)^T$ if $d=2$ and $\mathbf{x}=(x,y,z)^T$ if $d=3$, $\nabla$ is the gradient operator with respect to the spatial coordinate $\mathbf{x}\in\mathbb{R}^d$, $V(\mathbf{x})\geq0$ is a real-valued external potential, $\beta\in\mathbb{R}$ is a real parameter describing the strength of the interaction between particles (positive for repulsive interaction and negative for attractive interaction), $\Omega$ is an angular velocity, $\overline{\phi}$ denotes the complex conjugate of the wave function $\phi$, and
\begin{align*}
L_z=-i(x\partial_y-y\partial_x)
\end{align*}
is the $z$-component of the angular momentum operator $\mathbf{L}=\mathbf{x}\times\mathbf{P}$ with $\mathbf{P}=-i\nabla$ the momentum operator. We set $L_z=0$ and write $\mathbf{x}=x$ if $d=1$ in the model problem \eqref{eq:GPErot-energy}-\eqref{eq:ELeq-cons} to make it cover the (non-rotating) one-dimensional case.

This turns out to be a constrained non-convex minimization problem: Find $\phi_g\in S$ such that
\begin{align}\label{eq:gsdef}
E_g:=E(\phi_g)=\inf_{\phi\in S} E(\phi),
\end{align}
where $S=\left\{\phi\,|\, \|\phi\|^2=1,\,E(\phi)<\infty\right\}$ denotes the $L^2$ unit sphere. The associated Euler-Lagrange equation (or first-order optimality condition) is
\begin{align}\label{eq:ELeq}
\left\{\begin{aligned}
& \mu\phi(\mathbf{x}) = -\frac12\Delta\phi(\mathbf{x})+V(\mathbf{x})\phi(\mathbf{x})+\beta|\phi(\mathbf{x})|^2\phi(\mathbf{x})-\Omega L_z\phi(\mathbf{x}),\quad \mathbf{x}\in\mathbb{R}^d, \\
& \|\phi\|^2=1,
\end{aligned}\right.
\end{align}
which is a nonlinear eigenvalue problem for $(\mu,\phi)$ under the normalization condition, where $\Delta$ is the Laplace operator with respect to $\mathbf{x}$. The first equation in \eqref{eq:ELeq} is also known as the stationary/time-independent GP equation. Each eigenfunction $\phi$ of \eqref{eq:ELeq} is called a stationary state solution, and thus the ground state $\phi_g$ is the one with least energy. In addition, when $\phi\in S$ is an eigenfunction, the corresponding eigenvalue $\mu=\mu(\phi)\in\mathbb{R}$ is also called the chemical potential, which can be computed as
\begin{align}\label{eq:mu-chempot}
\mu(\phi)= \int_{\mathbb{R}^d} \left( \frac12|\nabla\phi|^2+V|\phi|^2+\beta|\phi|^4 -\Omega\overline{\phi}L_z\phi\right) \mathrm{d}\mathbf{x} = E(\phi) +\frac{\beta}{2}\int_{\mathbb{R}^d}|\phi(\mathbf{x})|^4\mathrm{d}\mathbf{x}.
\end{align}

Theoretical results on the existence of a minimizer to the minimization problem \eqref{eq:gsdef} can be found, e.g., in \cite{BC2013KRM,BWM2005CMS,LSY2000PRA,Seiringer2002CMP} and the references therein. It is noted that, for a non-rotating BEC, i.e. $\Omega=0$, the ground state is unique up to a constant phase and is strictly positive, while for a rotating BEC with suitably large rotation speed (e.g., $|\Omega|>\Omega^c$ for some critical rotation speed $\Omega^c>0$), problem \eqref{eq:gsdef} may no longer admits a unique minimizer \cite{BC2013KRM,BWM2005CMS}. On the other hand, existing experimental observations and numerical simulations \cite{ALT2017JCP,BWM2005CMS,CYRPB2018PRA,DP2017SISC,ECHC2002PRL,Cornell2003giant,ZZ2009CPC} show that very complicated vortex patterns (e.g., Abrikosov lattice and giant vortices) appear in the ground state of fast rotating BECs, especially in the strong repulsive interaction regime (i.e., $\beta\gg1$), and the corresponding GP energy landscape can be quite complex (e.g., it may have many local minima) and difficult to minimize. Therefore, one of the particular interests in numerical study of BECs is to develop numerical methods for computing the ground state of a rotating BEC with high computational efficiency and high accuracy.

Anchored by the above GP theory framework, in the past two decades, various feasible numerical methods have been proposed to compute the ground state of rotating or non-rotating BECs, mainly consisting of energy minimization methods using gradient flow \cite{BD2004SISC,BWM2005CMS,CST2000PRE,ZZ2009CPC,AD2014JCP,LC2021SISC,CL2021JCP,ZS2019JCP} or other optimization techniques \cite{ALT2017JCP,ATZ2018CiCP,BT2003JCP,DH2010JCP,GP2001SISC,DK2010SISC,DP2017SISC,HP2020SINUM,WWB2017JSC,CORT2009JCP,HSW2021JCP} and some nonlinear eigenvalue solvers \cite{Zhou2003NL,CCM2010SISC,WJC2013CiCP,AHP2021NM}. Here we mainly focus on energy minimization methods. In particular, we are going to briefly introduce the story of the normalized gradient flow methods since it is closely related to the work of this paper and is one of the most popular techniques for computing the ground states of BECs in the literature. In fact, the normalized gradient flow approaches, including the continuous normalized gradient flow (CNGF) and the gradient flow with discrete normalization (GFDN) or the imaginary time evolution method, was originally proposed for non-rotating case \cite{CST2000PRE,BD2004SISC,BCL2006JCP} and then extended to rotating case \cite{BWM2005CMS,ZZ2009CPC,AD2014JCP}. In 2005, Bao et al. \cite{BWM2005CMS} applied the CNGF method with a linearized backward Euler finite difference discretization (which can also be viewed as a full discretization of the GFDN) \cite{BD2004SISC} to compute the ground, symmetric and central vortex states for a rotating BEC. Then, in  \cite{ZZ2009CPC}, Zeng and Zhang adopted the GFDN method with the semi-implicit Euler Fourier pseudospectral discretization to simulate vortex lattice structures of condensate ground states in rapid rotating BECs. In \cite{AD2014JCP}, Antoine and Duboscq designed a preconditioned Krylov subspace iterative method to solve the linearized backward Euler Fourier pseudospectral scheme of the GFDN method for rotating BECs. Recently, the authors in \cite{LC2021SISC} revisited the GFDN for single- and multi-component BECs in non-rotational frame and pointed out that, except for some special cases (e.g., the linearized backward Euler scheme in single-component case), the temporal discretizations of the GFDN generally produces spurious ground states since the GFDN itself suffers from artificial splitting error terms that are of the same order of the time step size. They further proposed the normalized gradient flow with Lagrange multiplier (GFLM), which can be viewed as an approximation of the CNGF or a modified GFDN by introducing explicit Lagrange multiplier term in the gradient flow part, so that various temporal discretizations based on it can capture the correct ground state \cite{LC2021SISC,CL2021JCP}. In fact, the main idea in \cite{LC2021SISC,CL2021JCP} is quite general, and in particular, it is straightforwardly applicable to the rotational situations. In addition to gradient flow approaches, energy minimization methods have also been considered for computing the ground states of rotating BECs based on optimization techniques, such as the preconditioned nonlinear conjugate gradient method \cite{ALT2017JCP,ATZ2018CiCP,GTA2021CPC}, the Sobolev gradient method \cite{DK2010SISC}, the Riemannian conjugate gradient method \cite{DP2017SISC} and the regularized Newton method \cite{WWB2017JSC}.

The normalized gradient flow approaches mentioned above are considerably popular and widely used to compute the ground states of BECs due to the following benefits: (i) the form of the gradient flow is concise, which can be implemented easily via different discretization strategies; (ii) only the gradient information (or first-order variational derivative) of the objective functional $E(\phi)$ is needed, which makes the computation at each iteration easy to handle when a suitable discretization such as the explicit or semi-implicit Euler scheme is adopted. However, the limitation of the gradient flow methods is also quite explicit. Firstly, the slow convergence makes it sometimes very expensive in practice, especially for challenging optimization problems, such as the computation of ground states for rotating BECs with high rotation speed and/or strong interatomic interaction. Additionally, gradient flow stops whenever a stationary state (or critical point of the non-convex minimization problem) is reached, thus it may fail to converge to the ground state of rotating BECs. With these in mind, we are very much motivated to find methods that can retain the benefits of the normalized gradient flow for a rotating BEC but overcome the shortcomings to a certain extent.

We entrust the potential of second-order inertial dynamical flows, which is a recent topic initialized by the work of Su et al. \cite{SBC2016JMLR}. There a second-order ordinary differential equation (ODE) was proposed and its connection to Nesterov's accelerated gradient method \cite{Nesterov1983} was revealed. The acceleration property of such second-order flows were proven theoretically. Since then, there has been plenty of research focused on inertial dynamics by second-order ODEs, see for instance \cite{attouch2018fast,ABC2021JEMS,BDES2021FCM,ZhaHof20} and the references therein. Second-order ODE formulations have some advantages in comparison to the discrete acceleration algorithms (heavy ball method \cite{Pol64} and Nesterov's accelerated gradient method \cite{Nesterov1983}). First, it has been shown that the ODEs can recover the discrete gradient methods if certain discretization is applied. On the other hand, mechanism of the dynamical flows is more intuitive which makes the analysis of many properties of the acceleration algorithms more direct using the continuous ODE formulations. In particular, there are different ways to discretize the ODEs, which would generate variants of numerical algorithms, and some of them can have better numerical stability and efficiency. Among several directions of investigations, a very recent and comprehensive study of second-order flows as ODEs from theoretical point of view can be found in \cite{ABC2021JEMS}. In Section \ref{sec:pre}, we give a short review on the development of the inertial dynamics with respect to ODEs for convex optimization. Second-order flows have been considered recently also in the setting of partial differential equations (PDEs) \cite{benyamin2020accelerated,DHZ2021SIIMS}, which brings new challenges to the topic, especially to their theoretical understanding and numerical analysis. In \cite{DHZ2021SIIMS}, a class of geometric second-order quasilinear hyperbolic PDEs with applications in imaging are studied, and some of their analytical properties are proven, e.g., well-posedness, asymptotic behavior of solutions. Numerical analysis is still quite open in \cite{DHZ2021SIIMS}, and also in general to the topic of second-order nonlinear hyperbolic PDEs with dissipation. Having said that, studies in the literature mostly focus on convex optimization. In many applications, however, problems emerged often have to be formulated as non-convex optimization, i.e., the objective function is non-convex or the feasible set is non-convex, such as the ground state problem of rotating BECs \eqref{eq:gsdef} studied in this paper.

\subsection{Contribution and structure of the paper}
We propose two types of dissipative second-order nonlinear hyperbolic PDEs (systems), and they both are generalized inertial dynamics aiming for the constrained non-convex minimization problem \eqref{eq:gsdef}.
The first one intends to have the whole trajectory of the evolutionary PDE to satisfy the normalization constraint. To achieve this, a Lagrange multiplier term is introduced into the second-order PDE, which is capable of preserving the normalization constraint. This is a similar formulation as the CNGF \cite{BD2004SISC,BWM2005CMS}. Under this guideline, we use again an explicit update of the Lagrange multiplier coupled with a discrete normalization, as in the GFLM \cite{LC2021SISC}.
The second one is based on a different idea which we do not enforce the exact constraint for the whole trajectory but only at the end of convergence. This is realized by leveraging the formulation of augmented Lagrange multiplier for generating the dynamical flow, from which a coupled differential system with primal and dual variables is derived. These novel second-order hyperbolic PDEs (systems) become our starting point to develop efficient numerical algorithms. We choose an explicit scheme and also a semi-implicit scheme for the temporal discretization, in combination with a Fourier pseudospectral method for spatial discretization. Particularly, finite difference schemes of second-order accuracy for the temporal discretization (of first-order and second-order time derivatives) are employed to better simulate their continuous counterparts.  At every time step, either a simple update is needed for the explicit scheme, or a linear system with constant coefficients needs to be solved for the semi-implicit scheme, where some stabilization technique is applied. This also shows a difference between the continuous flows and discrete gradient methods, as a different discretization scheme would produce a variant numerical algorithm. {In the end, they contribute a rich family of much more efficient and robust algorithms in comparison with algorithms from gradient flow type methods: Larger sizes of time steps can be applied for explicit schemes, and a much smaller number of iterations are consumed for semi-implicit algorithm using the same step size. Note that for a fixed step size, the computational complexity of second-order flow algorithms at every time step is the same as the (first-order) gradient flow ones.} Moreover, the second-order flow based algorithms are observed to be less sensitive to the initial guess, and can reach lower energy state solutions, in particular, for fast rotating BECs. Plenty of numerical examples are provided for verification and comparison. In all cases, the methods corresponding to the second-order flows with suitable damping coefficients can perform significantly better than the methods corresponding to the normalized gradient flow. In particular, we test a few examples with strong interaction and fast rotation where the gradient flow approaches often consume a much larger unit of computational power and the Newton-type methods in the literature generally converge to some higher energy state solution, whereas the proposed methods produce the most competitive results.

During the preparation of our manuscript, we are aware that a damped second-order approach named dynamical functional particle method (DFPM) was proposed in \cite{OG2020EJP,GOJPA}, { which was based on an earlier work in \cite{EdvSveGulPer12}}. These works considered constrained nonlinear Schr\"odinger equations and applied to a rotating BEC with a different setting from \eqref{eq:gsdef}, in particular different dynamics and also different numerical methods were studied there. A further remark on the differences between DFPM and our proposed approaches can be found in Example~\ref{eg:5.4} in Section~\ref{sec:num}.

The rest of this paper is organized as follows. First, some preliminaries for second-order flows are provided in Section~\ref{sec:pre}. Then, two types of second-order flows with damping for the energy minimization problem \eqref{eq:gsdef} are introduced in Section~\ref{sec:proposed_flows}. In Section~\ref{sec:disc}, several efficient temporal and spatial discretization strategies are proposed to yield practical numerical algorithms. Detailed numerical results are reported in Section~\ref{sec:num} to illustrate the efficiency and accuracy of our algorithms. Finally, some concluding remarks are given in Section~\ref{sec:concl}.

We emphasize here that the concept of second-order flows should not be confused with second-order methods in optimization. The former means an evolutionary equation involves second-order derivative with respect to the (artificial) time variable, while the latter often requires second-order derivative (or its approximations) of the objective function in optimization, such as Newton's method. It should also be distinguished from the second-order accuracy for the numerical schemes in temporal discretization.

\section{Primer of second-order flows for convex optimization}
\label{sec:pre}

\begin{figure}[H]
	\subfigure[Steepest descent method]{
		\begin{minipage}[t]{0.52\linewidth}
			\centering
			\includegraphics[width=2.4in]{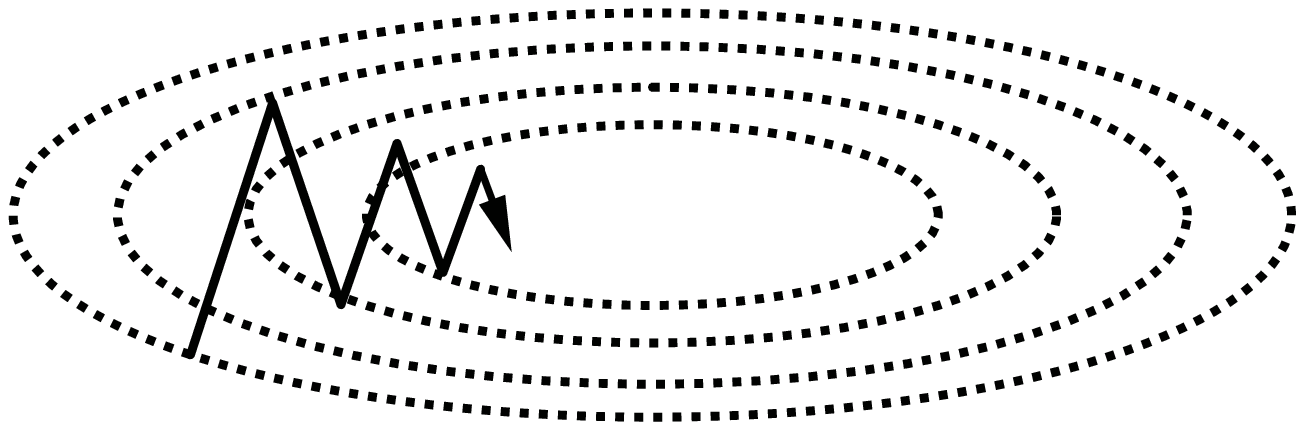}
			\label{fig1a}
	\end{minipage}}%
	\hspace{3mm}
	\subfigure[Heavy ball method]{
		\begin{minipage}[t]{0.42\linewidth}
			\centering
			\includegraphics[width=2.4in]{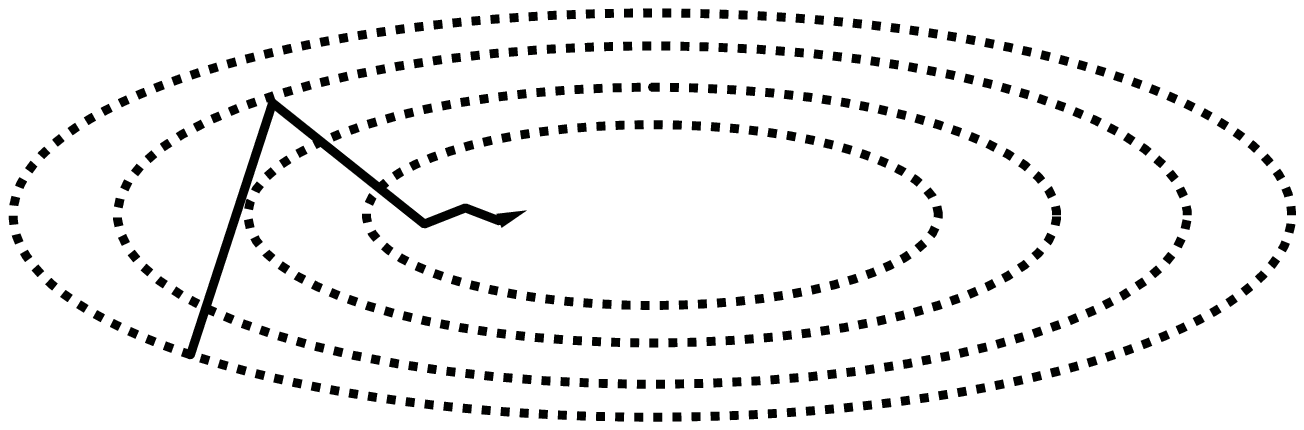}
			\label{fig1b}
	\end{minipage}}
	\caption{Illustration of the heavy ball method to alleviate the zig-zag phenomenon.}
	\label{fig1}
\end{figure}

In order to motivate our second-order flows for computing the ground state of a rotating BEC, here we briefly introduce some preliminaries of inertial dynamics in convex optimization. Let us pick the following abstract unconstrained minimization problem for illustration:
\begin{equation}\label{eq:opt_prob}
\min_{u\in H}\, f(u),
\end{equation}
where $f:H\to \mathbb{R}$ is a given convex differentiable function on a Hilbert space $H$. Our discussion starts with the steepest descent method (SDM) for minimizing  the convex function $f$, which reads
\begin{align}\label{SDM}
u_{k}=u_{k-1}-s \nabla f\left(u_{k-1}\right),\quad k=1,2,\ldots, 
\end{align}
with $s>0$ a step length and $u_0\in H$ a given initial guess. However, the SDM is not an optimal convergent first-order method, and the convergence rate can be very slow. In order to improve the convergence rates of the SDM, there has been a lot of research in the literature. One of them is the heavy ball method proposed by Polyak in 1964 \cite{Pol64}, which has the form
\begin{align}\label{HB}
u_{k+1}=u_{k}-s \nabla f\left(u_{k}\right)+\gamma\left(u_{k}-u_{k-1}\right),\quad k=1,2,\ldots,
\end{align}
for given $u_0$ and $u_1$. Here the term $\gamma\left(u_{k}-u_{k-1}\right)$, called a momentum, corrects the local steepest descent direction by extrapolating the direction with previous steps, and $\gamma>0$ is some constant. As shown in Fig.~\ref{fig1}, the negative gradient direction does not point towards the minimizer in most of the iterations. The SDM tends to oscillate from one side to the other, progressing slowly to the minimizer (see Fig.~\ref{fig1a}). Fig.~\ref{fig1b} shows how the momentum term helps to speed up convergence to the minimizer by damping these oscillations, even though theoretical optimal convergence rates seems not available for heavy ball method.
Nesterov's accelerated gradient method \cite{Nesterov1983} is based on a similar idea as heavy ball method, which uses, instead, an adaptive coefficient in front of the momentum term. It gives the following updates, and in particular returns improved convergence rates (optimal as first-order methods):
\begin{equation}\label{Nesterov}
\left\{\begin{aligned}
& u_{k}=w_{k-1}-s \nabla f\left(w_{k-1}\right), \\
& w_{k}=u_{k}+\frac{k-1}{k+2}\left(u_{k}-u_{k-1}\right).
\end{aligned}\right.
\end{equation}
These iterative type methods are fundamental in optimization, and more interestingly they can connect to some continuous evolutionary differential equations, which is summarized in Table~\ref{tab:first-order}.
\begin{table}[H]
\renewcommand{\arraystretch}{1.5}
	\centering
	\caption{Convergence rates of representative continuous flows and corresponding first-order optimization methods to the minimizer $u^{\star}$ of a convex function $f$.}
	\label{tab:first-order}
	\vspace{2mm}
	\small
	\resizebox{.95\textwidth}{!}{%
		\begin{tabular}{|l|l|l|l|}
			\hline
		Continuous flow & $f(u(t))-f(u^{\star})$ & Corresponding optimization method & $f(u_k)-f(u^{\star})$  \\ \hline \hline
		 $\dot{u}(t)=-\nabla f(u(t))$ & $O(1/t)$&	Steepest descent method \eqref{SDM}  &$O(1/k)$ \\ \hline
		 $\ddot{u}(t)+\eta \dot{u}(t)+\nabla f(u(t))=0$,   & 
		 \multirow{2}{*}{$O(1/t)$} & 	\multirow{2}{*}{Heavy ball method \eqref{HB}} & \multirow{2}{*}{$O(1/k)$} \\ 
		 where $\eta>0$ is a constant &&& \\ \hline
		 $\ddot{u}(t)+\frac{3}{t} \dot{u}(t)+\nabla f(u(t))=0$ & $O(1/t^{2})$ \cite{SBC2016JMLR} &	Nesterov's acceleration method \eqref{Nesterov}&  $O(1/k^{2})$ \cite{Nesterov1983} \\ \hline 
		\end{tabular}%
	}
	
\end{table}
In fact, in the literature, second-order flows with general damping coefficients $\frac{\alpha}{t}$ for $\alpha>0$ are often considered, e.g. \cite{SBC2016JMLR,attouch2018fast}. It was shown that $\alpha=3$ is a critical point, as for $\alpha<3$, the optimal convergence rate of first-order methods may not be possible in that case. For $\alpha>3$, improved convergence rates are available given some strong convexity assumptions on $f$.

It is not hard to find that the SDM can be viewed as an explicit Euler method for the gradient flow (i.e., the first continuous flow in Table~\ref{tab:first-order}).  Similarly, the heavy ball method and Nesterov's acceleration method are also corresponding to explicit time discretization schemes of the second \cite{attouch2000heavy} and the third continuous flows \cite{SBC2016JMLR}, respectively, in Table~\ref{tab:first-order}, which can be included into a common form as
\begin{align}\label{eq:SF}
\left\{\begin{aligned}
& \ddot{u}(t)+\eta(t)\dot{u}(t)=-\nabla f(u(t)), \quad t>0,\\
& u(0)=u_0,\quad \dot{u}(0)=v_0.
\end{aligned}\right.
\end{align}
Here $\eta(t)>0$ is the given damping coefficient, and $u_0$ and $v_0$ are the initial position and velocity, respectively. Note that $\eta(t)$ being a positive constant corresponds to the heavy ball method, and $\eta(t)=3/t$ corresponds to Nesterov's acceleration method. The results in Table~\ref{tab:first-order} have been much enriched or improved in the literature. Some recent comprehensive investigation within convex optimization can be found for instance in \cite{attouch2018fast,ABC2021JEMS}. Moreover, different aspects of second-order flows as regularization methods for linear ill-posed problems have been developed in \cite{BDES2021FCM}, where a detailed comparison of convergence rates of first-order and second-order flows are provided.

\begin{figure}[H]
	\centering	\includegraphics[width=0.38\textwidth]{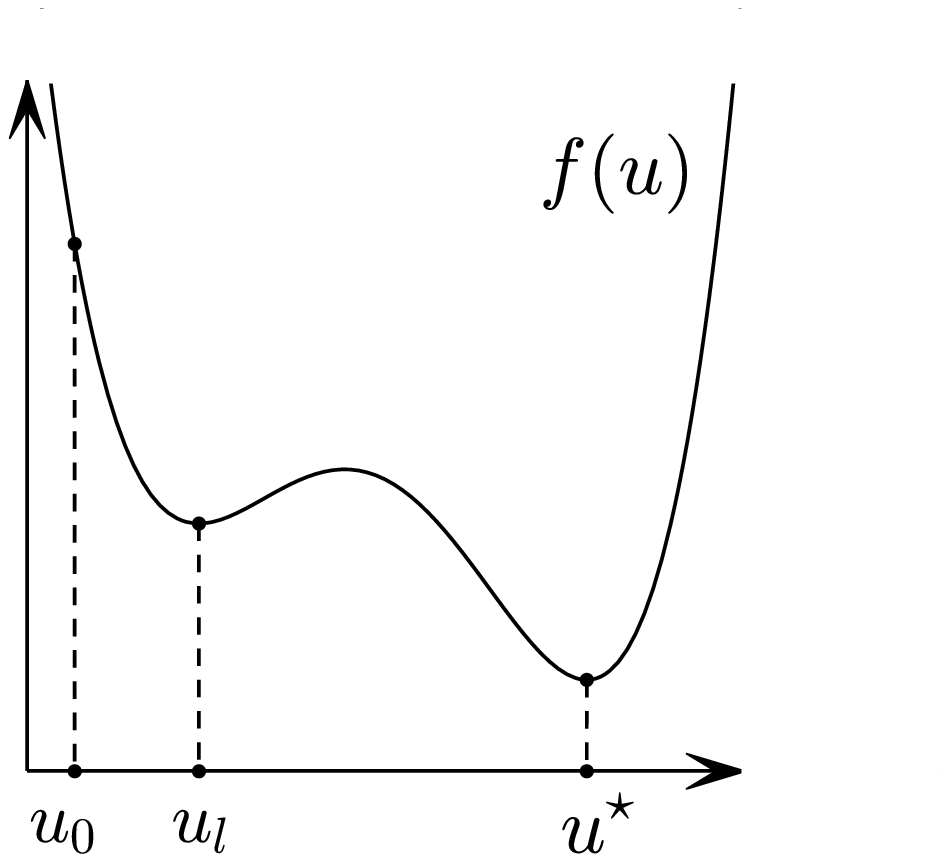}
	\vspace{-3ex}
	\caption{Profile of a one-dimensional non-convex function $f$ with local minima.}
	\label{fig:diagram}
\end{figure}

Even though the second-order flows can accelerate the convergence of gradient flows in some sense, we mention that they are in general not monotonically descent methods for minimizing the convex function $f$. Instead, the ODE \eqref{eq:SF} decays the ``total energy" or ``Hamiltonian" defined as $\mathcal{E}(t)=\frac{1}{2}\|\dot{u}(t)\|_H^2+f(u(t))$, of which the first term is similar to the ``kinetic energy" and the second one is viewed as the ``potential energy". At this point, we see that the second-order flows are actually fundamental objectives in classical mechanics. Another interesting difference between the gradient flow and the second-order flow is that the latter is capable of escaping saddle points or even local minima while the former typically stops at those points, when both are applied to non-convex optimization. This can be observed from the one-dimensional non-convex example shown in Fig.~\ref{fig:diagram}, where the non-convex objective function $f$ has a global minimizier $u^{\star}$ and a local minimizer $u_l$ with $f(u_l)>f(u^{\star})$. For the gradient flow of $f$ with the initial value $u_0$, the flow will typically approaching $u=u_l$, but for the second-order flow with proper choice of the initial velocity and the damping coefficient, it has the potential to reach the global minimizer $u=u^{\star}$ instead of $u=u_l$. 

The acceleration property and the potential of escaping local minima and/or saddle points of the second-order flow inspired our work in this paper, which aims to find solutions of the constrained non-convex minimization problem in \eqref{eq:gsdef}.

\section{Two families of damped second-order flows for computing BEC ground states}\label{sec:proposed_flows}
In this section, we propose two different second-order flows with respect to two different ideas of dealing with the normalization constraint. The first one uses a similar idea as the CNGF \cite{BD2004SISC,BWM2005CMS} and the GFLM \cite{LC2021SISC,CL2021JCP}, where we try to enforce the constraint to be satisfied along the whole trajectory of the second-order flow by embedding a Lagrange multiplier into the flow. The second one is built on the augmented Lagrangian with a primal-dual coupled system, and it does not require the constraint to be satisfied all the time but only when it comes to the end of the convergence. Due to the connection to normalized gradient flows, the conventional normalized gradient flow methods and the GFLM method to compute ground states of rotating BECs are described in \ref{sec:first-order} for comparison purposes. Throughout this paper, we use $\langle\cdot,\cdot\rangle$ to denote the $L^2$ inner product, i.e.,
\[ \langle u,v\rangle=\int_{\mathbb{R}^d} \overline{u(\mathbf{x})}v(\mathbf{x})\mathrm{d}\mathbf{x}, \quad\forall u,v\in L^2(\mathbb{R}^d), \]
and the $L^2$ norm $\|\cdot\|$ is reduced from this inner product. For simplicity, we assume $\beta\geq0$ in subsequent discussions.

\subsection{Damped second-order flow with discrete normalization}
Inspired by the second-order ODEs for convex optimization introduced in Section~\ref{sec:pre} and the CNGF method described in \eqref{eq:cngf}, we first propose a {\em continuous normalized second-order flow} for $\phi=\phi(\mathbf{x},t)$ to approach the ground state of a rotating BEC defined in \eqref{eq:gsdef} as
\begin{subequations}\label{eq:cndsf}
\begin{numcases}{}
\ddot{\phi} + \eta(t)\dot{\phi} 
=\frac12\Delta\phi-V\phi-\beta|\phi|^2\phi+\Omega L_z\phi+\lambda_{\phi}(t)\phi,
\quad \mathbf{x}\in\mathbb{R}^d,\; t>0, \label{eq:cndsf-a} \\
\phi(\mathbf{x},0)=\phi_0(\mathbf{x}),\quad \dot{\phi}(\mathbf{x},0)=v_0(\mathbf{x}), \quad \mathbf{x}\in\mathbb{R}^d,
\end{numcases}
\end{subequations}
where $t$ is the artificial time variable, $\eta(t)>0$ is a given damping coefficient, $\phi_0\in S$ is an initial guess for the ground state, and $v_0$ is the initial velocity. Two typical examples of the damping coefficient are $\eta >0$ being a constant and $\eta(t)=\alpha/t$ for some $\alpha\geq3$, which are respectively correspond to the heavy ball method and (generalized) Nesterov's acceleration method in convex optimization described in Section~\ref{sec:pre}.
The Lagrange multiplier $\lambda_{\phi}(t)$ intends to make the flow \eqref{eq:cndsf} preserve automatically the normalization constraint \eqref{eq:ELeq-cons}, i.e., $\|\phi(\cdot,t)\|^2\equiv\|\phi_0\|^2=1$ for all $t>0$. There are probably several formulations of $\lambda_{\phi}(t)$ to achieve this. Here we construct an explicit expression of $\lambda_{\phi}(t)$ according to the following two conditions that are expected to hold for a suitably smooth solution $\phi=\phi(\mathbf{x},t)$ of \eqref{eq:cndsf}:
\begin{align}\label{eq:lambda-cond1}
\frac{\mathrm{d}}{\mathrm{d}t}\left(\|\phi\|^2\right)
=2\mathrm{Re}\langle\phi,\dot{\phi}\rangle=0 \quad\mbox{and}\quad
\frac{\mathrm{d}^2}{\mathrm{d}t^2}\left(\|\phi\|^2\right)
= 2\left(\mathrm{Re}\,\langle\phi,\ddot{\phi}\rangle+\|\dot{\phi}\|^2\right)
=0,\quad \forall t>0. 
\end{align}
Taking the $L^2$ inner products with $\phi$ on both sides of \eqref{eq:cndsf-a} and then taking the real parts, yields
\begin{align}\label{eq:lambda-cond2}
\mathrm{Re}\langle \phi,\ddot{\phi}\rangle+\eta(t)\mathrm{Re}\langle \phi,\dot{\phi}\rangle
&=-\int_{\mathbb{R}^d}\left(\frac12|\nabla\phi|^2+V|\phi|^2+\beta|\phi|^4-\Omega\overline{\phi}L_z\phi\right)\mathrm{d}\mathbf{x}+\lambda_{\phi}(t)\|\phi\|^2.
\end{align}
Substituting \eqref{eq:lambda-cond1} into \eqref{eq:lambda-cond2}, we derive (formally)
\begin{align}\label{eq:lambda}
\lambda_{\phi}(t) 
&=\frac{\int_{\mathbb{R}^d}\left(\frac12|\nabla\phi|^2+V|\phi|^2+\beta|\phi|^4-\Omega\overline{\phi}L_z\phi\right)\mathrm{d}\mathbf{x}-\|\dot{\phi}(\cdot,t)\|^2}{\|\phi(\cdot,t)\|^2}. 
\end{align}
It is easy to see that the Lagrange multiplier $\lambda_{\phi}(t)$ given in \eqref{eq:lambda} is associated with the GP energy functional $E(\phi)$ \eqref{eq:GPErot-energy} and the Lagrange multiplier $\mu_{\phi}(t)$ in the CNGF \eqref{eq:cngf} through the following expression:
\begin{align}
\lambda_{\phi}(t)
=\frac{E(\phi)+\frac{\beta}{2}\int_{\mathbb{R}^d}|\phi(\mathbf{x},t)|^4\mathrm{d}\mathbf{x}-\|\dot{\phi}(\cdot,t)\|^2}{\|\phi(\cdot,t)\|^2}
=\mu_{\phi}(t)-\frac{\|\dot{\phi}(\cdot,t)\|^2}{\|\phi(\cdot,t)\|^2},
\end{align}
where $\mu_{\phi}(t)=\frac{1}{\|\phi(\cdot,t)\|^2}\int_{\mathbb{R}^d}\left(\frac12|\nabla\phi|^2+V|\phi|^2+\beta|\phi|^4-\Omega\overline{\phi}L_z\phi\right)\mathrm{d}\mathbf{x}$ is the same as in \eqref{eq:cngf}.
In addition, the condition \eqref{eq:lambda-cond1} also requires that the initial velocity $v_0$ satisfies $\mathrm{Re}\langle\phi_0,v_0\rangle=0$.

Formally, at the steady state or fixed point of the flow \eqref{eq:cndsf}, all the time derivative terms vanish, i.e., $\dot{\phi}=0$ and $\ddot{\phi}=0$, which leads to $\lambda_{\phi}(t)=\mu_{\phi}(t)$. With the normalization constraint $\|\phi\|^2=1$ conserved by the flow \eqref{eq:cndsf}, the steady state equation of \eqref{eq:cndsf-a} matches exactly the Euler-Lagrange equation \eqref{eq:ELeq}. Actually, the fixed point of \eqref{eq:cndsf} in the $(\phi,\dot{\phi})$-phase space is $(\phi,\dot{\phi})=(\phi_*,0)$ with $\phi_*$ an eigenfunction to the Euler-Lagrange equation \eqref{eq:ELeq}, i.e., a stationary state of the BEC.

We introduce the total cost function associated with the solution $\phi=\phi(\mathbf{x},t)$ to \eqref{eq:cndsf} defined as 
\begin{align}
\mathcal{F}(t) = E(\phi(\cdot,t))+\|\dot{\phi}(\cdot,t)\|^2.
\end{align}
The following property shows that, under certain conditions, the continuous normalized second-order flow \eqref{eq:cndsf} is indeed normalization-conservative and has a dissipation structure with respect to $\mathcal{F}(t)$.
\begin{property}
Let $\phi=\phi(\mathbf{x},t)$ be the solution to \eqref{eq:cndsf} with $\lambda_{\phi}(t)$ given in \eqref{eq:lambda} and initial data satisfying $\|\phi_0\|=1$ and $\mathrm{Re}\langle\phi_0,v_0\rangle=0$. {We further assume $\phi\in C^2\left((0,\infty), X\right)$, where $X=\left\{\phi \in H^1(\mathbb{R}^d):\, E(\phi)<\infty\right\}$.}
If the damping coefficient $\eta(\cdot)$ is positive and locally integrable on $(0,+\infty)$, 
then we have
\begin{align}
\|\phi(\cdot,t)\|^2\equiv\|\phi_0\|^2=1, \quad
\frac{\mathrm{d}}{\mathrm{d}t}\mathcal{F}(t) = -2\eta(t)\|\dot{\phi}(\cdot,t)\|^2,\quad
\forall\, t>0.
\end{align}
\end{property}
\begin{proof}
Combining \eqref{eq:lambda-cond2} and \eqref{eq:lambda}, we have $\frac{\mathrm{d}}{\mathrm{d}t}\mathrm{Re}\langle \phi,\dot{\phi}\rangle=-\eta(t)\mathrm{Re}\langle \phi,\dot{\phi}\rangle$. Since $\mathrm{Re}\langle \phi,\dot{\phi}\rangle\big|_{t=0}=\mathrm{Re}\langle\phi_0,v_0\rangle=0$ and $\eta>0$ is locally integrable on $(0,t]$ for all $t>0$, we conclude that
\[ \mathrm{Re}\langle\phi,\dot{\phi}\rangle=\lim\limits_{\varepsilon\rightarrow 0^{+}}\mathrm{e}^{-\int_{\varepsilon}^t\eta(s)\mathrm{d}s}\,\mathrm{Re}\langle\phi(\cdot,\varepsilon),\dot{\phi}(\cdot,\varepsilon)\rangle=0,\quad \forall t>0. \]
Therefore, $\frac{\mathrm{d}}{\mathrm{d}t}\left(\|\phi\|^2\right)=2\,\mathrm{Re}\langle\phi,\dot{\phi}\rangle=0$, and the conservation of the normalization constraint is verified. Further, utilizing \eqref{eq:cndsf-a} and noting that $\mathrm{Re}\langle \phi,\dot{\phi}\rangle=0$, a direct calculation shows that
\begin{align*}
\frac{\mathrm{d}}{\mathrm{d}t}\mathcal{F}(t) 
&= 2\,\mathrm{Re}\langle\dot{\phi},\ddot{\phi}\rangle + 2\,\mathrm{Re}\left\langle -\frac12\Delta\phi+V\phi+\beta|\phi|^2\phi-\Omega L_z\phi,\, \dot{\phi}\right\rangle \\
&= 2\,\mathrm{Re}\left\langle\dot{\phi},-\eta(t)\dot{\phi}+\lambda_{\phi}(t)\phi\right\rangle \\
&= - 2\,\eta(t)\|\dot{\phi}\|^2.
\end{align*}
Since $\eta(t)>0$, this gives the dissipation of $\mathcal{F}(t)$, and the proof is completed.
\end{proof}

From the decaying property of $\mathcal{F}(t)$ by the flow \eqref{eq:cndsf} described in above, we have (formally)
\begin{align}
E(\phi(\cdot,t))+\|\dot{\phi}(\cdot,t)\|^2+2\int_0^t\eta(s)\|\dot{\phi}(\cdot,s)\|^2\mathrm{d}s=\mathcal{F}(0)=E(\phi_0)+\|v_0\|^2,\quad \forall t>0.
\end{align}
Thus, the GP energy $E(\phi)$ along the trajectory of \eqref{eq:cndsf} is always bounded by $\mathcal{F}(0)=E(\phi_0)+\|v_0\|^2$. {In particular, if the initial velocity is taken as $v_0=0$, the whole trajectory of the flow \eqref{eq:cndsf} is naturally contained in the energy sublevel set $\{\phi\in S:\,E(\phi)\leq E(\phi_0)\}$. This provides some local stability to the method of the continuous normalized second-order flow \eqref{eq:cndsf} for computing the ground state. Note that when the initial value $\phi_0$ is taken as a metastable state, i.e., a local minimizer of the GP energy functional (if it exists), a nontrivial initial velocity $v_0$ may help to escape from this local minimum, which is another remarkable feature of the second-order flow methods compared to the ones by gradient flow.}

It is noted that the system \eqref{eq:cndsf} is a partial integro-differential equation for $\phi=\phi(\mathbf{x},t)$, in which the nonlocal Lagrange multiplier term $\lambda_{\phi}(t)\phi$ will bring difficulties for possible normalization-preserving discretizations. For the convenience of numerical implementation, referring to the idea of the GFLM \eqref{eq:gflm}, we propose the following computational model, which is referred to as the {\bfseries damped second-order flow with discrete normalization (DSFDN)}:
\begin{subequations}\label{eq:DSFDN}
\begin{numcases}{}
\ddot{\phi} + \eta(t)\dot{\phi} = \frac12\Delta\phi-V\phi-\beta|\phi|^2\phi+\Omega L_z\phi+\lambda_{\phi}(t_n)\phi,\quad t\in(t_n,t_{n+1}), \label{eq:DSFDN-a} \\
\phi(\mathbf{x},t_{n+1}):=\phi(\mathbf{x},t_{n+1}^+)=\frac{\phi(\mathbf{x},t_{n+1}^-)}{\|\phi(\cdot,t_{n+1}^-)\|},\quad n\geq0, \label{eq:DSFDN-b} \\
\dot{\phi}(\mathbf{x},t_{n+1}):=\dot{\phi}(\mathbf{x},t_{n+1}^+)=\dot{\phi}(\mathbf{x},t_{n+1}^-),\quad n\geq0, \label{eq:DSFDN-c}\\
\phi(\mathbf{x},0)=\phi_0(\mathbf{x}),\quad \dot{\phi}(\mathbf{x},0)=v_0(\mathbf{x}), \quad \mathbf{x}\in\mathbb{R}^d.
\end{numcases}
\end{subequations}
Here $t_n=n\tau$ ($n=0,1,\ldots$) and $\tau>0$ is a time step length, $\phi(\mathbf{x}, t_n^{\pm})=\lim_{t\to t_n^{\pm}} \phi(\mathbf{x},t)$, $\dot{\phi}(\mathbf{x}, t_n^{\pm})=\lim_{t\to t_n^{\pm}} \dot{\phi}(\mathbf{x},t)$, and $\lambda_{\phi}(t_n)$ is given explicitly as
\begin{align*}
\lambda_{\phi}(t_n)=\!\int_{\mathbb{R}^d}\!\left(\frac12|\nabla\phi(\cdot,t_n)|^2+V|\phi(\cdot,t_n)|^2+\beta|\phi(\cdot,t_n)|^4-\Omega\overline{\phi}(\cdot,t_n)L_z\phi(\cdot,t_n)\!\right)\!\mathrm{d}\mathbf{x}-\|\dot{\phi}(\cdot,t_n)\|^2.
\end{align*}

The system \eqref{eq:DSFDN}, as a piecewise-defined PDE system, can be viewed as an approximation of the continuous normalized second-order flow \eqref{eq:cndsf}. Formally, as $\tau\to0$, we have $\lambda_{\phi}(t)-\lambda_{\phi}(t_n)\to0$, then the second-order flow equation \eqref{eq:DSFDN-a} converges to \eqref{eq:cndsf} and the normalization factor \eqref{eq:DSFDN-b} converges to identity. Due to the fact that the evolution equation \eqref{eq:DSFDN-a} in each time interval is a hyperbolic PDE (with damping term), the initial conditions for both $\phi$ and $\dot{\phi}$ at $t=t_{n+1}$ ($n\geq0$) need to be specified whenever the evolution over a new time interval $(t_{n+1},t_{n+2})$ is considered. Various choices of the initial value of $\dot{\phi}$ at $t=t_{n+1}$ ($n\geq0$) are possible. As in \eqref{eq:DSFDN-c}, we simply use $\dot{\phi}(\cdot,t_{n+1}^{-})$, which is obtained by the evolution of \eqref{eq:DSFDN-a} in the previous time interval $(t_n,t_{n+1})$. 
An alternative way to define the initial condition of $\dot{\phi}$ at $t=t_{n+1}$ is to assign it as the projection of $\dot{\phi}(\cdot,t_{n+1}^{-})$ onto the $L^2$-orthogonal complement of $\phi(\cdot,t_{n+1})$. This makes the instantaneous rate of change of the normalization constraint at $t=t_{n+1}$ is zero (cf. \eqref{eq:lambda-cond1}) and seems to produce a good analog to the continuous flow \eqref{eq:cndsf} in terms of the normalization-preserving property. However, such an orthogonal projection is somewhat cumbersome or redundant in the DSFDN, because the normalization condition at $t=t_{n+1}$ is guaranteed strictly by \eqref{eq:DSFDN-b}. Thus, we only focus on the current version of \eqref{eq:DSFDN-c} for brevity and convenience. Actually, as we will see later, the condition \eqref{eq:DSFDN-c} does not have to appear explicitly in further discretizations of the DFSDN.

{We emphasize that the explicit Lagrange multiplier term $\lambda_{\phi}(t_n)\phi$ in \eqref{eq:DSFDN-a} plays a crucial role to ensure that the DSFDN \eqref{eq:DSFDN} can capture the correct stationary state with a finite time step length $\tau>0$. 
To illustrate this, we simply assume that $\phi(\cdot,0)=\phi_g$ and $\dot{\phi}(\cdot,0)=0$ in \eqref{eq:DSFDN}, where $\phi_g\in S$ denotes the ground state defined in \eqref{eq:gsdef}. It implies
\begin{align*}
\lambda_{\phi}(0)=
\int_{\mathbb{R}^d}\left(\frac12|\nabla\phi_g|^2+V|\phi_g|^2+\beta|\phi_g|^4-\Omega\,\overline{\phi}_g L_z\phi_g\right)\mathrm{d}\mathbf{x}
=\mu(\phi_g)=:\mu_g,
\end{align*}
the chemical potential of $\phi_g$ \eqref{eq:mu-chempot}, and the evolution equation \eqref{eq:DSFDN-a} becomes a steady dynamics at $(\phi,\dot{\phi})=(\phi_g,0)$ due to the fact that $(\mu_g,\phi_g)$ solves the Euler-Lagrange equation \eqref{eq:ELeq}. Further, the normalization factor in the projection step \eqref{eq:DSFDN-b} becomes $\|\phi(\cdot,t_{n+1}^-)\|=1$. On the other hand, as a twin version of \eqref{eq:DSFDN} we consider the situation without the explicit Lagrange multiplier term $\lambda_{\phi}(t_n)\phi$ in \eqref{eq:DSFDN-a}. Then the resulting evolutionary equation starting from the ground state $\phi_g$ with zero initial velocity $v_0=0$ leads to a non-steady dynamics at $(\phi,\dot{\phi})=(\phi_g,0)$, and due to the nonlinear effects from the interaction term $\beta|\phi|^2\phi$ (when $\beta\neq0$), it usually fails to pull $\phi(\cdot,t_{n+1}^-)$ back to $\phi_g$ with the subsequent normalization step that is a linear process. Thus, in general, dropping the explicit Lagrange multiplier term $\lambda_{\phi}(t_n)\phi$ in \eqref{eq:DSFDN-a} causes the system \eqref{eq:DSFDN} to produce an incorrect solution with an artificially introduced error depending on the time step parameter $\tau>0$. In addition, the explicit Lagrange multiplier term $\lambda_{\phi}(t_n)\phi$ in \eqref{eq:DSFDN-a} can be taken as other forms, e.g., $\lambda_{\phi}(t_n)\phi(\cdot,t_n)$, as long as $(\phi,\dot{\phi})=(\phi_g,0)$ is the fixed point of the system \eqref{eq:DSFDN}. }

In terms of dealing with constraint, the DSFDN method \eqref{eq:DSFDN} uses a normalization at each time step to exactly project the solution back to the constraint manifold $S$. For the minimization problem of a single-component BEC model with only one normalization constraint on which we focus in this paper, the exact projection is rather easy to handle. This might not be the case when one considers the ground state problem of general multi-component BECs with multiple constraints, e.g., high-spin BEC system with the two constraints on the total mass and magnetization, as the computation of the exact projection rises some challenges \cite{CL2021JCP}. In the next subsection, we propose another framework to compute the ground state of a rotating BEC, which can be better generalized in the face of more complex constraints.

\subsection{Damped second-order flow based on augmented Lagrangian method}
Augmented Lagrange multiplier is a very popular framework in constrained optimization, which combines the strengths of Lagrange multiplier and penalty method. The advantage is that it stabilizes the update of the primal variable using the augmented term but without sending the penalty parameter to $\infty$ to have the constraint to be satisfied. Based on the formulation of the augmented Lagrangian, we consider here another novel damped second-order dynamical systems. To adapt better to the normalization constraint of the non-convex minimization problem \eqref{eq:gsdef}, we propose the following general form of augmented Lagrangian 
\begin{align}\label{augmented Lagrangian}
\mathcal{L}_{\sigma}(\phi, \chi):=E(\phi)-\chi (\|\phi\|^2-1)+\sigma R(\|\phi\|^2-1),
\end{align}
where $\chi\in\mathbb{R}$ is the multiplier and the augmentation parameter $\sigma>0$ is a given constant. Here $R(\cdot):\mathbb{R}\to \mathbb{R}$ is some penalty function which satisfies
\begin{enumerate}[(i)]
  \item $R(0)=0$ and $R(\zeta)\geq 0$ for all $\zeta\in \mathbb{R}$, and it is coercive in the positive direction, i.e., $\lim_{\zeta\to+\infty} \frac{R(\zeta)}{\zeta}=+\infty$.
  \item $R$ is convex, continuously differentiable, and $R'(0)=0$. 
\end{enumerate}
In this paper, we will consider two examples of such functions:
\begin{equation}
(a) \; R(\zeta)=\frac{1}{2}|\zeta|^2 \qquad \text{ and } \qquad 
(b) \; R(\zeta)=\frac{1}{2}(\max(\zeta,0))^2. 
\end{equation}
Standard augmented Lagrangian usually takes the quadratic penalty function in Example $(a)$. In the normalization constraint case, the quadratic penalty appears to be not the best choice, as the functional $(\|\phi\|^2-1)^2$ is non-convex, which makes the optimization problem harder to solve in particular when the penalty parameter $\sigma$ gets larger and larger. To make the penalty term to be convex, Example $(b)$ seems to be more appealing for the normalization constraint.
In our numerical experiments in Section~\ref{sec:num}, we employ Example $(b)$ for two-dimensional rotating BECs.

The Lagrangian dual function is then defined as
\begin{align}\label{eq:dual-fun}
g(\chi)=\inf \limits_{\phi}\mathcal{L}_{\sigma}(\phi, \chi),
\end{align}
and its corresponding dual problem is
\begin{align}\label{eq:dual-pb}
\sup \limits_{\chi}g(\chi).
\end{align}
This converts the primal problem \eqref{eq:gsdef} to an {\bfseries unconstrained} saddle-point problem for $(\phi,\chi)$, i.e.,
\begin{align}
\sup_{\chi} \inf_{\phi} \mathcal{L}_{\sigma}(\phi, \chi),
\end{align}
and its optimal value gives a lower bound on the optimal value of the primal problem.

For brevity, hereafter we use
\begin{align}\label{eq:DE}
\frac{\delta E(\phi)}{\delta\overline{\phi}}=-\frac12\Delta\phi+V\phi+\beta|\phi|^2\phi-\Omega L_z\phi
\end{align}
to denote the Wirtinger's variational derivative of the GP energy functional $E(\phi)$ \eqref{eq:GPErot-energy} with respect to a complex-valued function $\phi:\mathbb{R}^d\to\mathbb{C}$. Let $(\phi_{\star}, \chi_{\star})$ be a saddle point of the augmented Lagrangian $\mathcal{L}_{\sigma}(\phi, \chi)$. Then $(\phi_{\star}, \chi_{\star})$ satisfies the Karush-Kuhn-Tucker (KKT) condition
\begin{align}\label{eq:AL-kkt}
\left\{\begin{aligned}
& \frac{\delta \mathcal{L}_{\sigma}}{\delta\overline{\phi}}(\phi_{\star}, \chi_{\star})
=\frac{\delta E}{\delta\overline{\phi}}(\phi_{\star})-\chi_{\star}\phi_{\star}+\sigma R^\prime(\|\phi_{\star}\|^2-1)\phi_{\star}
=0,\\
& \frac{\partial\mathcal{L}_{\sigma}}{\partial\chi}(\phi_{\star}, \chi_{\star})=\|\phi_{\star}\|^{2}-1=0,
\end{aligned}\right.
\end{align}
which is exactly the nonlinear eigenvalue problem \eqref{eq:ELeq} for $(\chi_{\star}, \phi_{\star})$ since $R'(0)=0$. 

To approach the saddle point of the augmented Lagrangian $\mathcal{L}_{\sigma}(\phi, \chi)$, we first consider the following gradient flow based on augmented Lagrangian (GFAL) for $(\phi, \chi)=(\phi(\mathbf{x},t), \chi(t))$:
\begin{subequations}\label{gfal}
\begin{numcases}{}
\dot{\phi}=-\frac{\delta E}{\delta\overline{\phi}}(\phi) +\chi(t)\phi-\sigma R^\prime(\|\phi\|^2-1)\phi, \quad \mathbf{x}\in\mathbb{R}^d,\; t>0, \label{eq:gfal-a}\\
\frac{1}{\xi(t)}\dot{\chi} = 1-\|\phi\|^{2}, \quad  t>0,\label{eq:gfal-b} \\
\phi(\cdot,0)=\phi_0,\quad \dot{\phi}(\cdot,0)=v_0,\quad \chi(0)=\chi_0,
\end{numcases}
\end{subequations}
where $\xi(t)>0$ is a real function which serves as a scaling factor to control the pace of the evolution of $\chi(t)$, and $(\phi_0, v_0, \chi_0)$ are the initial data. Intuitively, the negative gradient of the penalty (i.e., $-\sigma R^\prime(\|\phi\|^2-1)\phi$ in \eqref{eq:gfal-a})  behaves as a pull-back force that makes $\phi$ move towards the normalization constraint manifold $S$ with a magnitude proportional to the augmentation parameter $\sigma$. When $\|\phi\|^2\gg 1$, due to the coercive property of $R$, this pull-back force can become significant (in fact, it grows superlinearly as $\|\phi\|\to\infty$). When the normalization constraint is approaching, $\dot{\chi}$ will also gradually go to $0$.

Slightly different to the inertial dynamics of gradient flows for the primal problem \eqref{eq:gsdef}, we introduce the following primal-dual system called {\bfseries damped second-order flow based on augmented Lagrangian (DSFAL)}: 
\begin{subequations}\label{eq:dsfal}
\begin{numcases}{}
\ddot{\phi}+\eta(t)\dot{\phi}=-\frac{\delta E}{\delta\overline{\phi}}(\phi) +\chi(t)\phi-\sigma R^\prime(\|\phi\|^2-1)\phi, \quad \mathbf{x}\in\mathbb{R}^d,\; t>0, \label{eq:dsfal-a}\\
\frac{1}{\xi(t)}\dot{\chi} = 1-\|\phi\|^{2}, \quad t>0, \label{eq:dsfal-b} \\
\phi(\cdot,0)=\phi_0,\quad \dot{\phi}(\cdot,0)=v_0,\quad \chi(0)=\chi_0,
\end{numcases}
\end{subequations}
with a given damping coefficient $\eta(t)>0$. {Clearly, when the system comes to a steady state, i.e., $\dot{\phi}$, $\dot{\chi}$ and $\ddot{\phi}$ all become 0, $(\phi, \chi)$ satisfies exactly the KKT condition \eqref{eq:AL-kkt}.} Note that it couples a second-order flow for primal variable $\phi$ and a first-order flow for dual variable $\chi$. The heuristic of \eqref{eq:dsfal} is that an inertial dynamics of $\mathcal{L}_{\sigma}(\phi,\cdot)$ is applied with respect to $\phi$, whereas a gradient flow is still employed to evolve $\chi$. We notice that $\mathcal{L}_{\sigma}(\cdot,\chi)$ is formally a linear function with respect to $\chi$ whose gradient direction is the steepest rising direction. In \cite[Theorem 2.8]{DHZ2021SIIMS}, it is proven that there exists no finite extinction time of second-order flows for homogeneous functional, whereas the extinction time is finite for gradient flows of homogeneous functionals (though it was proven for a special case, the idea is identical for general homogeneous functionals). This shows that in the linear function case, first-order flows are superior in terms of convergence rates, which also gives the intuition why we choose first-order flow for the multiplier $\chi$. In our numerical experiments, we find that this second-order and first-order coupled system for $\phi$ and $\chi$ has better performance in comparing with the coupled system with both second-order flows for primal and dual variables. The latter has been recently proposed in the literature, e.g., \cite{he2021convergence}.
\begin{remark}
We notice that the second-order system \eqref{eq:dsfal} may require different scaling factor $\xi(t)$ from the first-order system \eqref{gfal}. For the latter, $\xi(t)$ can be chosen as a positive constant, while for the former, an attenuation function, which satisfies $\xi(t)\to 0$ as $t\to +\infty$ and $\int_{0}^\infty \xi(t)dt=+\infty$, is proposed to improve the numerical stability. Later in the numerical tests in Section \ref{sec:num}, we consider $\xi(t)$ of the following form for both first- and second-order systems \eqref{gfal}-\eqref{eq:dsfal}:
\begin{align}\label{xi-form}
   \xi(t)=\left\{\begin{aligned}
& b, & t\leq t_s, \\
& b/(t-t_s+1), & t>t_s.
\end{aligned}\right.
\end{align}
Here $b>0$ is a constant, and $t_s>0$ is a turning point which we choose to be dependent on the constraint violation $|\|\phi\|^2-1|$. Such a decay property of $\xi(t)$ appears to be crucial for the computational stability of the second-order flow methods while it does not matter much for the gradient flow methods. However, we do not investigate this topic theoretically here.
\end{remark}

A decaying total cost function for the DSFAL may be obtained under certain convexity assumptions from some existing works \cite{luo2021primal,he2021convergence}. However, for the non-convex problem considered in this paper, whether it is possible to construct a decaying cost function for the DSFAL is still open, which we do not pursue here. Nevertheless, our numerical results in Section~\ref{sec:num} demonstrate that the DSFAL is feasible and efficient for our constrained non-convex minimization problem \eqref{eq:gsdef}.

The two proposed PDEs (systems) \eqref{eq:cndsf} and \eqref{eq:dsfal} turn out to be novel, and they are interesting topics in PDE analysis on their own right. However, analytical aspects of the two PDEs, e.g., well-posedness and regularity of solutions, are out of the scope of the current paper, we leave them for future work and for readers who are interested in these problems.

\section{Discretization of the proposed second-order flows}
\label{sec:disc}
In this section, we provide several numerical strategies to discretize the two types of second-order flows proposed in the last section.
We investigate the temporal and spatial discretization strategies separately, with some emphases on the former. In terms of temporal discretization, we give an explicit discretization scheme (the leap-frog scheme) and a type of semi-implicit discretization method. Both schemes are of second-order accuracy in terms of approximation rates to time partial derivatives, which contribute to the stability of the corresponding algorithms. A Fourier pseudospectral method is applied in all the cases for spatial discretization.

In the subsequent discussions, we assume $v_0 \equiv 0$ in both DSFDN \eqref{eq:DSFDN} and DSFAL \eqref{eq:dsfal} for simplicity.

\subsection{Temporal discretizations of DSFDN}\label{sec:disc-DSFDN}
We provide two temporal discretization schemes for the DSFDN \eqref{eq:DSFDN}. One is an explicit leap-frog scheme and the other is a semi-implicit scheme with stabilization.

Set $\phi^0=\phi_0$ and $t_n=n\tau$ with $\tau>0$ a given time step length. Let $\phi^n$ be a numerical approximation of $\phi(\cdot,t_n)$ for $n=1,2,\ldots$. Denote $\eta^n=\eta(t_n)=\eta(n\tau)$, $n=1,2,\ldots$, and
\begin{align}\label{eq:Gphin}
G(\phi^{n})=\frac{\delta E}{\delta\overline{\phi}}( \phi^{n})=-\frac12\Delta{\phi}^{n}+V\phi^n+\beta|\phi^n|^2\phi^n-\Omega L_z\phi^n,\quad n=0,1,\ldots.
\end{align}
In the following temporal discretizations for the DSFDN \eqref{eq:DSFDN}, the explicit Lagrange multiplier will be approximated by
\begin{align}\label{eq:DFSDN-lambda-n}
\lambda^n:=\begin{cases}
\mu^0-\|v_0\|^2=\mu^0,& n=0, \\
\mu^n-\|(\phi^{n}-\phi^{n-1})/\tau\|^2, & n=1,2,\ldots,
\end{cases}
\end{align}
with
\begin{align*}
\mu^n: = \left\langle G(\phi^{n}),\phi^n\right\rangle=
\int_{\mathbb{R}^d}\left(\frac12|\nabla\phi^n|^2+V|\phi^n|^2+\beta|\phi^n|^4-\Omega\overline{\phi^n}L_z\phi^n\right)\mathrm{d}\mathbf{x}
=E(\phi^n)+\frac{\beta}{2}\int_{\mathbb{R}^d}|\phi^n|^4\mathrm{d}\mathbf{x}.
\end{align*}
Here the forward first-order finite difference is applied to approximate $\dot{\phi}(\cdot,t_n)$ in the Lagrange multiplier $\lambda_{\phi}(t_n)$ in \eqref{eq:DSFDN}. It is believed that there are many other reasonable ways to do this. Our numerical results presented in Section~\ref{sec:num} show that the current version of $\lambda^n$ is feasible and accurate. More precisely, for the two schemes described below, the error of the numerical solution $\phi^n$ of the DSFDN is indeed of second-order accuracy with respect to the time step size $\tau$. We will, then, stick to the form of $\lambda^n$ in \eqref{eq:DFSDN-lambda-n} through out the paper.

\subsubsection{Leap-frog temporal discretization} 
We start with the explicit leap-frog scheme for the DSFDN \eqref{eq:DSFDN} (DSFDN-LF) given by
\begin{subequations}\label{eq:nsflm-fe}
\begin{align}
& \frac{\tilde{\phi}^{n+1}-2\phi^n+\phi^{n-1}}{\tau^2}+\eta^n\frac{\tilde{\phi}^{n+1}-\phi^{n-1}}{2 \tau} = -G(\phi^n)+\lambda^n\phi^n,\quad \mathbf{x}\in\mathbb{R}^d,\\
& \phi^{n+1}= \tilde{\phi}^{n+1}/\|\tilde{\phi}^{n+1}\|,\quad n\geq 1.
\end{align}
\end{subequations}
Since this is a three-level scheme, to initialize, we compute the first step $\phi^1\approx \phi(\cdot,\tau)$ by using the second-order Taylor's expansion at $t=0$ and noting $\dot{\phi}(\cdot,0)=v_0=0$:
\begin{subequations}\label{eq:si-firststep}
\begin{align}
& \tilde{\phi}^1=\phi^0+\tau\dot{\phi}\big|_{t=0}+\frac{\tau^2}{2}\ddot{\phi}\big|_{t=0} =\phi^0+\frac{\tau^2}{2}\left(-G(\phi^0)+\lambda^0\phi^0\right),\\ 
& \phi^1=\tilde{\phi}^1/\|\tilde{\phi}^1\|.
\end{align}
\end{subequations}
In detail, the DSFDN-LF \eqref{eq:nsflm-fe}-\eqref{eq:si-firststep} can be summarized as in Algorithm~\ref{alg:1}.
\begin{algorithm}[!ht]
    \caption{Leap-frog scheme for DSFDN (DSFDN-LF)}
    \label{alg:1}
    \begin{algorithmic}
\STATE {\bf First step:} Compute $\tilde{\phi}^1=\phi^0+\frac{\tau^2}{2}\left(-G(\phi^0)+\lambda^0\phi^0\right)$ and $\phi^1=\tilde{\phi}^1/\|\tilde{\phi}^1\|$
        \WHILE {not converge}
            \STATE Compute $G(\phi^n)$ and $\lambda^n$ according to \eqref{eq:Gphin} and \eqref{eq:DFSDN-lambda-n}, respectively
            \STATE Compute $\tilde{\phi}^{n+1}$ via
            \begin{align}\label{eq:leap-frog}
            \tilde{\phi}^{n+1}=\frac{(\tau\eta^n-2)\phi^{n-1}+4\phi^{n}+2\tau^2(-G(\phi^n)+\lambda^n\phi^n)}{2+\tau\eta^n} 
            \end{align}
            \STATE Update
            $\phi^{n+1}= \tilde{\phi}^{n+1}/\|\tilde{\phi}^{n+1}\|$
            \STATE $n:=n+1$
        \ENDWHILE
    \end{algorithmic}
\end{algorithm}

\begin{remark}
Notice that \eqref{eq:leap-frog} can be rewritten as
\begin{align*}
\tilde{\phi}^{n+1}=\phi^n-\frac{2\tau^2}{2+\tau\eta^n}\Big(G(\phi^n)-\lambda^n\phi^n\Big)+\frac{2-\tau\eta^n}{2+\tau\eta^n}\Big(\phi^n-\phi^{n-1}\Big).
\end{align*}
This formulation reminds us again the connection of the second-order flows to the Nesterov's (or Polyak's) methods, which gives an interpretation of the accelerating phenomenon of second-order flows. {For convex optimization problems, in comparison with the gradient flow, the second-order flow method allows larger step sizes in explicit discretization, i.e., $2/\ell_h$ for gradient flows can be improved to $2/\sqrt{\ell_h}$ for the second-order flows, where $\ell_h$ is a relatively large constant determined by the spatial discretization and the Lipschitz constant of nonlinearities (see, e.g., \cite{benyamin2020accelerated}). For the non-convex optimization problem studied in this paper, we think that similar properties should hold locally. In our numerical experiments later, in term of computational stability, it is evident that the DSFDN-LF \eqref{eq:nsflm-fe} allows larger time step compared to the explicit scheme of gradient flow type methods such as the GFLM-FE \eqref{eq:gflm-fe}, though further rigorous mathematical discussions are not provided here.}
\end{remark}

The following property states that the convergent state of the DSFDN-LF scheme \eqref{eq:nsflm-fe}-\eqref{eq:si-firststep} is exactly the eigenfunction to the Euler-Lagrange equation \eqref{eq:ELeq}.

\begin{property}
The iterative scheme \eqref{eq:nsflm-fe} reaches its fixed point (i.e., $\phi^{n+1}=\phi^{n}=\phi^{n-1}$) for some $n\geq1$ if and only if $G(\phi^n)=\lambda^n\phi^n$ (or equivalently, $G(\phi^n)=\mu^n\phi^n$ and $\phi^{n}=\phi^{n-1}$).
\end{property}
\begin{proof}
{\em Necessity.} 
Suppose $\phi^{n+1}=\phi^{n}=\phi^{n-1}$ for some $n\geq1$. Then $\lambda^n=\mu^n$ and $\tilde{\phi}^{n+1}=c\,\phi^n$ with $c=\|\tilde{\phi}^{n+1}\|>0$, so that the iterative scheme \eqref{eq:nsflm-fe} becomes
\begin{align}\label{eq:nsflm-fe-limit}
\frac{(c-1)\phi^n}{\tau^2}+\eta^n\frac{(c-1)\phi^n}{2\tau}=-G(\phi^n)+\mu^n\phi^n.
\end{align}
Taking the $L^2$ inner product in both side of \eqref{eq:nsflm-fe-limit} with $\phi^n$ and noting that $\|\phi^n\|=1$, we get
\[ (c-1)\left(\frac{1}{\tau^2}+\eta^n\frac{1}{2\tau}\right)=-\langle G(\phi^n),\phi^n\rangle+\mu^n=0, \]
which implies $c=1$. Applying \eqref{eq:nsflm-fe-limit} and the fact $\lambda^n=\mu^n$, we arrive at $G(\phi^n)=\lambda^n\phi^n$.

{\em Sufficiency.} 
Suppose $G(\phi^n)=\lambda^n\phi^n$. Then, by utilizing $\|\phi^n\|=1$, we have $\lambda^n=\langle G(\phi^n),\phi^n\rangle=\mu^n$. Applying \eqref{eq:DFSDN-lambda-n} leads to $\phi^n=\phi^{n-1}$. Thus the iterative scheme \eqref{eq:nsflm-fe} becomes
\[
\frac{\tilde{\phi}^{n+1}-\phi^n}{\tau^2}+\eta^n\frac{\tilde{\phi}^{n+1}-\phi^n}{2\tau}=0.
\]
This means that $\tilde{\phi}^{n+1}=\phi^n$, and therefore, $\phi^{n+1}=\tilde{\phi}^{n+1}=\phi^n=\phi^{n-1}$.
\end{proof}

Since the proof of the above property is based on inner products, it is easy to check that the same result holds for a Galerkin-type full discretization of the DSFDN-LF scheme \eqref{eq:nsflm-fe}-\eqref{eq:si-firststep}.

According to the above property and noting that $\|\phi^n\|=1$ by the normalization step, the condition $G(\phi^n)=\mu^n\phi^n$ (i.e., $(\mu^n,\phi^n)$ solving the Euler-Lagrange equation \eqref{eq:ELeq}) is necessary, but not sufficient, for the DSFDN-LF scheme \eqref{eq:nsflm-fe}-\eqref{eq:si-firststep} to converge. This actually shows an essential difference between the second-order flow and gradient flow methods. In our numerical experiments, we terminate the DSFDN-LF iteration if the maximal residual of the Euler–Lagrange equation \eqref{eq:ELeq} at $(\mu^n,\phi^n)$,
\begin{align}\label{eq:residual-cre}
&e_{r}^{n}:= 
\left\|G( \phi^{n})-\mu^n \phi^{n}\right\|_{\infty}<\varepsilon_{r}
\end{align}
and the discretized velocity
\begin{align}\label{eq:velocity-cre}
&d_{v}^{n}:= 
\frac{\left\|\phi^{n}-\phi^{n-1}\right\|_{\infty}}{\tau}<\varepsilon_{v}
\end{align}
with $\varepsilon_{r}$ and $\varepsilon_{v}$ two tolerances. As shown in the above property, one can also compute the difference between three consecutive iteration points (e.g., $\max\{\|\phi^{n+1}-\phi^n\|,\|\phi^n-\phi^{n-1}\|\}$) or the residual of the Euler–Lagrange equation \eqref{eq:ELeq} at $(\lambda^n,\phi^n)$ (i.e., $\|G(\phi^n)-\lambda^n\phi^n\|_{\infty}$) to check the convergence of the DSFDN-LF scheme \eqref{eq:nsflm-fe}-\eqref{eq:si-firststep}.

\subsubsection{Stabilized semi-implicit temporal discretization}
We then propose to apply a semi-implicit scheme with stabilization for \eqref{eq:DSFDN} that reads as
\begin{subequations}\label{eq:CN-LF}
\begin{align}
&\frac{\tilde{\phi}^{n+1}-2\phi^n+\phi^{n-1}}{\tau^2}+\eta^n\frac{\tilde{\phi}^{n+1}-\phi^{n-1}}{2\tau}  \nonumber \\
&\qquad = \left(\frac12\Delta-\vartheta^n\right)\frac{\tilde{\phi}^{n+1}+\phi^{n-1}}{2}+\vartheta^n\phi^n-V\phi^n-\beta|\phi^n|^2\phi^n+\Omega L_z\phi^n+\lambda^n\phi^n,\\
& \phi^{n+1}= \tilde{\phi}^{n+1}/\|\tilde{\phi}^{n+1}\|, \quad n\geq1,
\end{align}
\end{subequations}
where $\vartheta^n\geq0$ is a stabilization factor. Similar to the stabilized semi-implicit schemes for normalized gradient flow approaches (cf. \cite{BC2013KRM,LC2021SISC}), the scheme is expected to be able to take a larger stable time step through a proper choice of the stabilization parameter $\vartheta^n$. {The same stabilization parameter as in the gradient flow case is adopted in our numerical experiments (see \eqref{eq:stabn} for the detail) and is shown to be effective numerically, though an optimal choice of it needs in-depth investigation, which could be an interesting topic for further study.} The first step $\phi^1$ in \eqref{eq:CN-LF} is also computed through \eqref{eq:si-firststep}, the same as in the DSFDN-LF scheme. 

Similar to the previous property of the DSFDN-LF scheme, we have the following result. The proof is omitted here for brevity.
\begin{property}
The iterative scheme \eqref{eq:CN-LF} reaches its fixed point (i.e., $\phi^{n+1}=\phi^{n}=\phi^{n-1}$) for some $n\geq1$ if and only if $G(\phi^n)=\lambda^n\phi^n$ (or equivalently, $G(\phi^n)=\mu^n\phi^n$ and $\phi^{n}=\phi^{n-1}$).
\end{property}
%


In practice, the DSFDN-SI scheme \eqref{eq:CN-LF} can be implemented as in Algorithm~\ref{alg:2} with the stopping criterion \eqref{eq:residual-cre}-\eqref{eq:velocity-cre}.

\begin{algorithm}[!ht]
    \caption{Semi-implicit scheme for DSFDN (DSFDN-SI)}
    \label{alg:2}
    \begin{algorithmic}
\STATE {\bf First step:} Compute $\tilde{\phi}^1=\phi^0+\frac{\tau^2}{2}\left(-G(\phi^0)+\lambda^0\phi^0\right)$ and $\phi^1=\tilde{\phi}^1/\|\tilde{\phi}^1\|$
        \WHILE {not converge}
        	\STATE Compute $G(\phi^n)$ and $\lambda^n$ according to \eqref{eq:Gphin} and \eqref{eq:DFSDN-lambda-n}, respectively
            \STATE Compute
        	\begin{align*}
             \mathcal{H}^n =\!\left(-1+\frac{\tau}{2}\eta^n-\frac{\tau^2}{2}\vartheta^n\!\right)\!\phi^{n-1}+\frac{\tau^2}{4}\Delta\phi^{n-1} +2\phi^n+\tau^2\left(\!(\vartheta^n+\lambda^n)\phi^n-\frac{1}{2}\Delta\phi^n-G(\phi^n)\!\right)
            \end{align*}
        	\STATE Solve the linear elliptic equation for $\tilde{\phi}^{n+1}$:
            \begin{align}\label{eq:linear-sys1}
            \left(1+\frac{\tau}{2}\eta^n+\frac{\tau^2}{2}\vartheta^n\right) \tilde{\phi}^{n+1}-\frac{\tau^2}{4}\Delta\tilde{\phi}^{n+1} =\mathcal{H}^n.
            \end{align}
            \STATE Update $\phi^{n+1}= \tilde{\phi}^{n+1}/\|\tilde{\phi}^{n+1}\|$
            \STATE $n:=n+1$
        \ENDWHILE
    \end{algorithmic}
\end{algorithm}

\begin{remark}
The main computational cost at each time step for the DSFDN-SI scheme is to solve \eqref{eq:linear-sys1}, a linear elliptic equation with constant coefficient, which can be solved efficiently by some fast Poisson solver, e.g., the fast Fourier transform (FFT), of which the details is presented in Section~\ref{sec:disc-spatial}. In fact, both DSFDN-SI and DFSDN-LF have $O(\tau^2)$ accuracy in temporal approximation when they are stable for the selected time step size $\tau$. This is confirmed by our numerical experiments (see Fig.~\ref{fig1.b} in Section~\ref{sec:num}). On the other hand, for the DSFDN-SI, the choice of time step size is less restrictive on the mesh size of the spatial discretization, thus obtaining a more stable algorithm without much computational cost increment compared to the DSFDN-LF when a fast Poisson solver is employed. To keep track with the trajectories of the PDEs, we choose the time step size $\tau$ to not larger than $0.1$ in most of our numerical experiments. 
\end{remark}

\subsection{Temporal discretizations of DSFAL}\label{sec:disc-DSFAL}
Taking into account the stability, approximation accuracy and computational efficiency, we focus on a similar semi-implicit discretization scheme as \eqref{eq:CN-LF} for DSFAL \eqref{eq:dsfal}, referred as DSFAL-SI:
\begin{subequations}\label{eq:dsfal-si}
\begin{align}
&\frac{{\phi}^{n+1}-2\phi^n+\phi^{n-1}}{\tau^2}+\eta^n\frac{{\phi}^{n+1}-\phi^{n-1}}{2\tau}=\left(\frac12\Delta-\vartheta^n\right)\frac{\phi^{n+1}+\phi^{n-1}}{2} \nonumber \\
&\qquad\qquad+\vartheta^n\phi^n-V\phi^n-\beta|\phi^n|^2\phi^n+\Omega L_z\phi^n+\chi^n\phi^n-\sigma R^\prime(\|\phi^n\|^2-1)\phi^n, \label{eq:dsfal-si-a}\\
&\frac{\chi^{n+1}-\chi^{n}}{\xi^n \tau}= 1-\|\phi^{n+1}\|^2,\quad n=1,2,\ldots, \label{eq:dsfal-si-b}
\end{align}
\end{subequations}
with $\vartheta^n\geq0$ a stabilization factor. $\xi^n$ evaluates from $\xi(t^n)$ using the form \eqref{xi-form}, and $t_s$ in \eqref{xi-form} is taken as $t_s=n_s \tau$, where $n_s$ represents the smallest positive integer $n$ satisfying
\begin{align} \label{xi-con}
    |\|\phi^n\|^2-1|<\epsilon_s,
\end{align}
where $\epsilon_s$ is some threshold value of small scales. Later, in the numerical examples in Section \ref{sec:num}, we take $\epsilon_s$ as $10^{-14}$. Starting from the initial data $\phi^0=\phi_0$ and $\chi^0=\chi_0$, we compute $\phi^1$ based on the second-order Taylor's expansion at $t=0$ and the initial velocity $\dot{\phi}(\cdot,0)=v_0=0$ as
\begin{align*}
\phi^1=\phi^0+\tau\dot{\phi}\big|_{t=0}+\frac{\tau^2}{2}\ddot{\phi}\big|_{t=0} =\phi^0+\frac{\tau^2}{2}\left(-G(\phi^0)+\chi^0\phi^0\right),
\end{align*}
and then set $\chi^{1}=\chi^{0}+\tau\xi^0(1-\|\phi^{1}\|^2)$. 
The specific iteration steps of the DSFAL-SI scheme \eqref{eq:dsfal-si} are summarized in Algorithm \ref{alg:3}.
\begin{algorithm}[!ht]
    \caption{Semi-implicit scheme for DSFAL (DSFAL-SI)}
    \label{alg:3}
    \begin{algorithmic}
    \STATE {\bf First step:} Compute $\phi^1=\phi^0+\frac{\tau^2}{2}\left(-G(\phi^0)+\chi^0\phi^0\right)$ and $\chi^{1}=\chi^{0}+\tau\xi^0(1-\|\phi^{1}\|^2)$
    \WHILE {not converge}
    	\STATE Compute
        	\begin{align*}
             \mathcal{H}^n &=\left(-1+\frac{\tau}{2}\eta^n-\frac{\tau^2}{2}\vartheta^n\right)\phi^{n-1}+\frac{\tau^2}{4}\Delta\phi^{n-1} \\ &\quad\;+2\phi^n+\tau^2\left(\Big(\vartheta^n+\chi^n-\sigma R^\prime(\|\phi^n\|^2-1)\Big)\phi^n-\frac{1}{2}\Delta\phi^n-G(\phi^n)\right)
            \end{align*}
        \STATE Update $\phi^{n+1}$ by solving the linear elliptic equation
            \begin{align*}
            \left(1+\frac{\tau}{2}\eta^n+\frac{\tau^2}{2}\vartheta^n\right)\phi^{n+1}-\frac{\tau^2}{4}\Delta\phi^{n+1} =\mathcal{H}^n.
            \end{align*}
        \STATE Update $\chi^{n+1}=\chi^{n}+\tau\xi^n(1-\|\phi^{n+1}\|^2)$ 
        \STATE $n:=n+1$
    \ENDWHILE
    \end{algorithmic}
\end{algorithm}

\begin{property}
The iterative scheme \eqref{eq:dsfal-si} reaches its fixed point (i.e., $\phi^{n+1}=\phi^{n}=\phi^{n-1}$ and $\chi^{n+1}=\chi^n$) for some $n\geq1$ if and only if the following hold:
\[ G(\phi^n)=\chi^n\phi^n,\quad \|\phi^n\|^2=1,\quad \phi^n=\phi^{n-1}. \]
Moreover, in this case, $\chi^n=\langle G(\phi^n),\phi^n\rangle=\mu^n$.
\end{property}
\begin{proof}
{\em Necessity.} 
Suppose $\phi^{n+1}=\phi^{n}=\phi^{n-1}$ and $\chi^{n+1}=\chi^n$ for some $n\geq1$. From \eqref{eq:dsfal-si-b}, we have $\|\phi^{n+1}\|^2=1$. Then $\|\phi^{n}\|^2=1$ and \eqref{eq:dsfal-si-a} becomes exactly $G(\phi^n)=\chi^n\phi^n$.

{\em Sufficiency.} 
Suppose $G(\phi^n)=\chi^n\phi^n$, $\|\phi^n\|^2=1$ and $\phi^n=\phi^{n-1}$ for some $n\geq1$. Then \eqref{eq:dsfal-si-a} becomes
\[
\left(\frac{1}{\tau^2}+\frac{\eta^n}{2\tau}+\frac{\vartheta^n}{2}-\frac14\Delta\right)\Big({\phi}^{n+1}-\phi^n\Big)=0,
\]
which implies $\phi^{n+1}=\phi^{n}$. Thus $\|\phi^{n+1}\|=\|\phi^n\|=1$, and \eqref{eq:dsfal-si-b} shows $\chi^{n+1}=\chi^{n}$.
\end{proof}

In our numerical experiments, the stopping conditions for the DSFAL-SI \eqref{eq:dsfal-si} is set as 
\begin{subequations}\label{eq:stop-pd}
\begin{align}
&e_{r}^{n}:= 
\left\|G( \phi^{n})-\mu^{n} \phi^{n}\right\|_{\infty}<\varepsilon_{r}, \\
&d_{v}^{n}:= 
\frac{\left\|\phi^{n}-\phi^{n-1}\right\|_{\infty}}{\tau}<\varepsilon_{v}, \\
&e_{c}^{n}:=|\|\phi^n\|^2-1|<\varepsilon_{c},
\end{align}
\end{subequations}
with some tolerances $\varepsilon_r$, $\varepsilon_v$ and $\varepsilon_c$, and the initial value for the Lagrange multiplier $\chi^0$ is taken as
\begin{align}
&\chi^0=\chi_0:=
\frac{\int_{\mathbb{R}^d}\left(\frac{1}{2}|\nabla\phi_0|^2+V|\phi_0|^2+\beta|\phi_0|^4-\Omega\overline{\phi}_0L_z\phi_0\right)\mathrm{d}\mathbf{x}}{\|\phi_0\|^2}+\sigma R^\prime(\|\phi_0\|^2-1).
\end{align}

For comparison, we describe the following semi-implicit scheme (GFAL-SI) to discretize the GFAL \eqref{gfal}:
\begin{subequations}\label{eq:gfal-si}
\begin{align}
\frac{\phi^{n+1}-\phi^{n}}{\tau}&=
\left(\frac12\Delta-\vartheta^n\right)\phi^{n+1}
+\bigg(\vartheta^n-V-\beta|\phi^n|^2+\Omega L_z+\chi^n-\sigma R^\prime(\|\phi^n\|^2-1)\bigg)\phi^n, \\
\frac{\chi^{n+1}-\chi^{n}}{\xi^n\tau}&= 1-\|\phi^{n+1}\|^2,\quad n=0,1,\ldots
\end{align}
\end{subequations}

\subsection{Spatial discretization with Fourier pseudospectral method}\label{sec:disc-spatial}
There are various options to further spatially discretize the previously proposed temporal semi-discrete schemes to produce corresponding full-discretization algorithms, such as finite element methods, finite difference methods or spectral methods. We adopt a Fourier pseudospectral method \cite{BC2013KRM,ZZ2009CPC} here. For simplicity, the following description is presented for two-dimensional cases ($d=2$), and the cases of other dimensions can be straightforwardly generalized.

Due to the external trapping potential $V(\mathbf{x})$, the stationary state solution decays to zero exponentially when $|\mathbf{x}|\rightarrow\infty$ \cite{BC2013KRM}. Therefore, in practical computations, the problem in $\mathbb{R}^d$ can be naturally truncated to a suitably large bounded domain, which is referred to as the computational domain here. Especially, we choose $\mathcal{D}=[-L, L]^2$ to be the computational domain with periodic boundary conditions, where $L>0$ is large enough so that the domain truncation error can be ignored. Uniformly sampled grid is adopted with mesh size $h=2L/M$, where $M$ is an even number. Then the domain $\mathcal{D}$ is divided with $(M+1)\times(M+1)$ grid points $\{(x_j,y_k)\}_{j,k=0}^M$ with $x_{j}=-L+j h$ and $y_{k}=-L+k h$ ($j,k=0,1,\ldots,M$).

Let $\phi_{jk}$ be the numerical approximation for the value of a function $\phi$ at $(x_j,y_k)$. The Fourier pseudospectral approximations to the operators $\Delta$ and $L_z$ are given as
\begin{align*}
&(\Delta_h \phi)_{jk}=-\sum_{p=-M/2}^{M/2 -1} \sum_{q=-M / 2}^{M/2-1}\left((\varrho^{x}_{p})^{2}+(\varrho^{y}_{q})^{2}\right) \widehat{\phi}_{pq}\, e^{i \frac{2j p \pi}{M}} e^{i \frac{2k q \pi}{M}}, \\
&(L_{z}^h\phi)_{jk}=-i\Big(x_{j} (\partial_{y}^{h}\phi)_{jk}-y_{k} (\partial_{x}^{h}\phi)_{jk}\Big), \qquad j,k=0,1,\cdots, M-1,
\end{align*}
where $\varrho^{x}_{p}=p \pi/L$, $\varrho^{y}_{q}=q \pi/L$, $-M/2\leq p,q\leq M/2-1$, and 
\begin{align*}
& (\partial_{x}^{h}\phi)_{jk}=\sum_{p=-M/2}^{M/2-1} \sum_{q=-M/2}^{M/2-1} i\varrho^{x}_{p} \widehat{\phi}_{pq}\, e^{i \frac{2j p \pi}{M}} e^{i \frac{2k q \pi}{M}}, \quad
(\partial_{y}^{h}\phi)_{jk}=\sum_{p=-M/2}^{M/2-1} \sum_{q=-M/2}^{M/2-1} i\varrho^{y}_{q} \widehat{\phi}_{pq}\, e^{i \frac{2j p \pi}{M}} e^{i \frac{2k q \pi}{M}},
\end{align*}
with $\widehat{\phi}_{pq} =\frac{1}{M^2} \sum_{j=0}^{M-1} \sum_{k=0}^{M-1} \phi_{jk}\, e^{-i \frac{2j p \pi}{M}} e^{-i \frac{2k q \pi}{M}}$ the discrete Fourier coefficients of mesh function $\phi_{jk}$. Therefore, the discretized energy $E_{h}(\cdot)$, chemical potential $\mu_{h}(\cdot)$ and $L^2$-norm $\|\cdot\|_h$ can be, respectively, computed as
\begin{align*}
& E_{h}(\phi)=h^2 \sum_{j=0}^{M-1} \sum_{k=0}^{M-1} \left(\frac12|(\partial_{x}^{h}\phi)_{jk}|^2+\frac12|(\partial_{x}^{h}\phi)_{jk}|^2+V(x_j,y_k)\left|\phi_{j k}\right|^{2}+\frac{\beta}{2}|\phi_{jk}|^{4}-\Omega\,\overline{\phi}_{jk}(L^{h}_z\phi)_{jk}\!\right)\!, \\
& \mu_h(\phi)=
E_{h}(\phi)+\frac{\beta h^2}{2} \sum_{j=0}^{M-1} \sum_{k=0}^{M-1}\left|\phi_{j k}\right|^{4}, \qquad
\|\phi\|_h^{2}=h^2 \sum_{j=0}^{M-1} \sum_{k=0}^{M-1}\left|\phi_{j k}\right|^{2}.
\end{align*}

Finally, we shall solve the linear elliptic equation with constant coefficient using the discrete Fourier transform (DFT), which is the main cost for the implementation of semi-implicit temporal discretization schemes presented in Section~\ref{sec:disc-DSFDN} and Section~\ref{sec:disc-DSFAL}. Taking \eqref{eq:linear-sys1} as an example to apply the Fourier pseudospectral spatial discretization with $d=2$, we get 
\begin{align}\label{eq:linear-sys-discrete}
    \left(1+\frac{\tau}{2}\eta^n+\frac{\tau^2}{2}\vartheta^n\right)\tilde{\phi}^{n+1}_{jk} - \frac{\tau^2}{4}\left(\Delta_h\tilde{\phi}^{n+1}\right)_{jk}= \mathcal{H}^n_{jk},\quad j,k=0,1,\ldots,M-1,
\end{align}
where $\tilde{\phi}^{n+1}_{jk}$ and $\mathcal{H}^n_{jk}$ are, respectively, the discretized values of $\tilde{\phi}^{n+1}$ and $\mathcal{H}^n$ in \eqref{eq:linear-sys1}. Then, performing the DFT on both sides of \eqref{eq:linear-sys-discrete}, yields
\begin{align*}
     \left(1+\frac{\tau}{2}\eta^n+\frac{\tau^2}{2}\vartheta^n+\frac{\tau^2}{4}\Big((\varrho^{x}_{p})^2+(\varrho^{y}_{q})^2\Big)\right)\widehat{(\tilde{\phi}^{n+1})}_{pq}= \widehat{(\mathcal{H}^n)}_{pq},\quad p,q=-\frac{M}{2},\ldots,\frac{M}{2}-1.
 \end{align*}
Taking the inverse DFT, the system is solved as
 \begin{align*}
\tilde{\phi}_{jk}^{n+1}=\sum_{p=-M/2}^{M/2 -1} \sum_{q=-M / 2}^{M/2-1}\frac{\widehat{(\mathcal{H}^n)}_{pq}}{1+\frac{\tau}{2}\eta^n+\frac{\tau^2}{2}\vartheta^n+\frac{\tau^2}{4}\left((\varrho^{x}_{p})^2+(\varrho^{y}_{q})^2\right)}\, e^{i \frac{2j p \pi}{M}} e^{i \frac{2k q \pi}{M}}.
 \end{align*}
We can see that the memory cost is $O(M^2)$ and the computational cost is $O(M^2\ln(M^2))$ when the FFT method is applied here.

\begin{remark} \label{multigrid}
Direct numerical computation with a very fine spatial grid in 2D and 3D is expensive, even when the FFT method is applied. This typically arises, for example, in the simulation of the ground state of fast rotating BECs, where a relatively small spatial mesh size is taken to capture complex vortex patterns with high resolution. Similar to \cite{WWB2017JSC}, we adopt a cascaded multigrid technique to further improve the computational efficiency. In 2D, the basic idea of its implementation is as follows: First, we solve the minimization problem \eqref{eq:gsdef} by an iterative algorithm presented in Section~\ref{sec:disc-DSFDN}-\ref{sec:disc-DSFAL} on a coarse spatial grid with $(2^p+1)\times(2^p+1)$ points for some positive integer $p$; Then we interpolate the obtained numerical ground state solution to the refined grid with $(2^{p+1}+1)\times(2^{p+1}+1)$ points to get the initial guess for the iterative algorithms on this finer grid; Repeat this process several times until we obtain the numerical ground state solution on the finest grid. This multigrid technique is applied in our numerical experiments of Example~\ref{eg:5.6} in Section~\ref{sec:num}.
\end{remark}

\section{Numerical results}\label{sec:num}
In this section, we show several numerical examples ranging from 1D to 2D to illustrate the accuracy and efficiency of the proposed methods. {Section 5.1 gives comparison of the different methods for 1D cases, while Section 5.2 is for 2D cases with and without the rotation term, and Section 5.3 sticks to 2D cases with rotation term under two different trapping potential functions.} All algorithms were implemented in Python (v3.8.10) and all experiments were performed on a workstation with a 2.40 GHz CPU.

We choose the damping coefficient $\eta(t)$ in the DSFDN \eqref{eq:DSFDN} and the DSFAL \eqref{eq:dsfal} as $\eta(t)=\alpha/t$, where $\alpha\geq3$ is a parameter that can be tuned for optimal numerical performances. {Initial velocity $v_0$ is uniformly taken as $v_0\equiv0$ throughout the experiments.} Three types of trapping potential $V(\mathbf{x})$ are involved in our numerical experiments: (i) the harmonic potential \cite{BC2013KRM}
\begin{align}\label{eq:harmonic-potential}
V(\mathbf{x}) = V_{\mathrm{ho}}(\mathbf{x}) :=
\frac{1}{2}\begin{cases}
\gamma_{x}^{2} x^{2}, & d=1, \\
\gamma_{x}^{2} x^{2}+\gamma_{y}^{2} y^{2}, & d=2, \\
\end{cases}
\end{align}
where $\gamma_x, \gamma_y>0$ are the trapping frequencies; (ii) the harmonic-plus-lattice potential \cite{ALT2017JCP,BC2013KRM}
\begin{align}\label{eq:harmonic-p-o-l}
V(\mathbf{x})=V_{\mathrm{ho}}(\mathbf{x})+\frac{\kappa}{2} \begin{cases}
\sin ^{2}\left(q_{x} x\right), & d=1,\\
\sin ^{2}\left(q_{x} x\right)+\sin ^{2}\left(q_{y} y\right), & d=2, \\
\end{cases}
\end{align}
where $\kappa$, $q_x$ and $q_y$ are positive constants; and (iii) the harmonic-plus-quartic potential in 2D \cite{WWB2017JSC}
\begin{align}\label{eq:harmonic-p-q}
V(\mathbf{x})=(1-\theta) \frac{x^2+y^2}{2}+\frac{\kappa\left(x^{2}+y^{2}\right)^{2}}{4}
\end{align}
with $\theta$ and $\kappa$ two positive constants. For non-rotating BECs, refer to \cite{ALT2017JCP}, we uniformly take the initial value $\phi_0$ as
\begin{align} \label{init_tf}
\phi_{0}(\mathbf{x})=\frac{\phi^{\mathrm{TF}}(\mathbf{x})}{\|\phi^{\mathrm{TF}}\|}\quad \text{with}   \quad  
\phi^{\mathrm{TF}}(\mathbf{x})=
\begin{cases}
\sqrt{\left(\mu^{\mathrm{TF}}-V(\mathbf{x})\right)/\beta}, & \mbox{if } V(\mathbf{x}) \leq \mu^{\mathrm{TF}}, \\ 0, & \mbox{otherwise},
\end{cases}
\end{align}
where $\mu^{\mathrm{TF}}=\frac{1}{2}\left(3 \beta \gamma_{x}\right)^{2/3}$ if $d=1$ and $\mu^{\mathrm{TF}}=\left(\beta \gamma_{x} \gamma_{y}\right)^{1/2}$ if $d=2$.

For fair comparisons, we conduct all numerical experiments on the bounded domain $\mathcal{D}$ with periodic boundary conditions, and the Fourier pseudospectral method described in Section~\ref{sec:disc-spatial} is utilized for spatial discretization with mesh size $h$. Furthermore, the stabilization parameters $\vartheta^n$ in all semi-implicit schemes are chosen as the following form:
\begin{align}\label{eq:stabn}
\vartheta^n=\frac12\Big((V+\beta|\phi^n|^2)_{\max}+(V+\beta|\phi^n|^2)_{\min}\Big),
\end{align}
where $(\cdot)_{\max}$ and $(\cdot)_{\min}$ represent respectively the maximum and minimum with respect to the spatial variable $\mathbf{x}\in\mathcal{D}$. 
Finally, the stopping criterion \eqref{eq:residual-cre}-\eqref{eq:velocity-cre} is adopted for the DSFDN-LF, DSFDN-SI, GFLM-FE and GFLM-BF schemes, while the stopping condition given in \eqref{eq:stop-pd} is applied for the DSFAL-SI and GFAL-SI schemes. In the following, we use ``Iter" and ``CPU(s)" to represent the number of iterations and the CPU time consumed to achieve convergence, respectively.

\subsection{Numerical results and comparisons in 1D examples}
In this subsection, we perform some numerical computations in the 1D situation to compare the numerical efficiency of the second-order flow methods proposed in Section~\ref{sec:proposed_flows}-\ref{sec:disc} with the GFLM methods described in \ref{sec:first-order}. The harmonic potential \eqref{eq:harmonic-potential} and the harmonic-plus-lattice potential \eqref{eq:harmonic-p-o-l} in 1D case are considered here with $\gamma_x = 1$, $\kappa = 25$ and $q_x = \frac{\pi}{2}$. $\mathcal{D}=[-32,32]$ is set to be the computational domain, and $h=\frac{1}{32}$ is taken as the mesh size. We choose the tolerances $\varepsilon_v=\varepsilon_r=10^{-12}$ and $\varepsilon_c=10^{-14}$ in \eqref{eq:residual-cre}-\eqref{eq:velocity-cre} and \eqref{eq:stop-pd} all through 1D cases. Unless specified, we use $E_g$ and $\mu_g$ to represent the corresponding GP energy and chemical potential of the computed ground state solution, respectively.

\begin{example}\label{5.1}
In this example, we investigate the properties of second-order flow algorithms through a simple example in 1D with harmonic potential \eqref{eq:harmonic-potential} and $\beta=250$. 
\end{example}
\vspace{-8pt}

\begin{figure}[!t]
\hspace*{2em}
\subfigure[The order of time error of DSFDN-SI with different $\alpha$ for approximating $\phi(\cdot,t)$ at $t=5$.]{
\begin{minipage}[t]{0.45\linewidth}
\centering
\includegraphics[width=2.9in,height=2.6in]{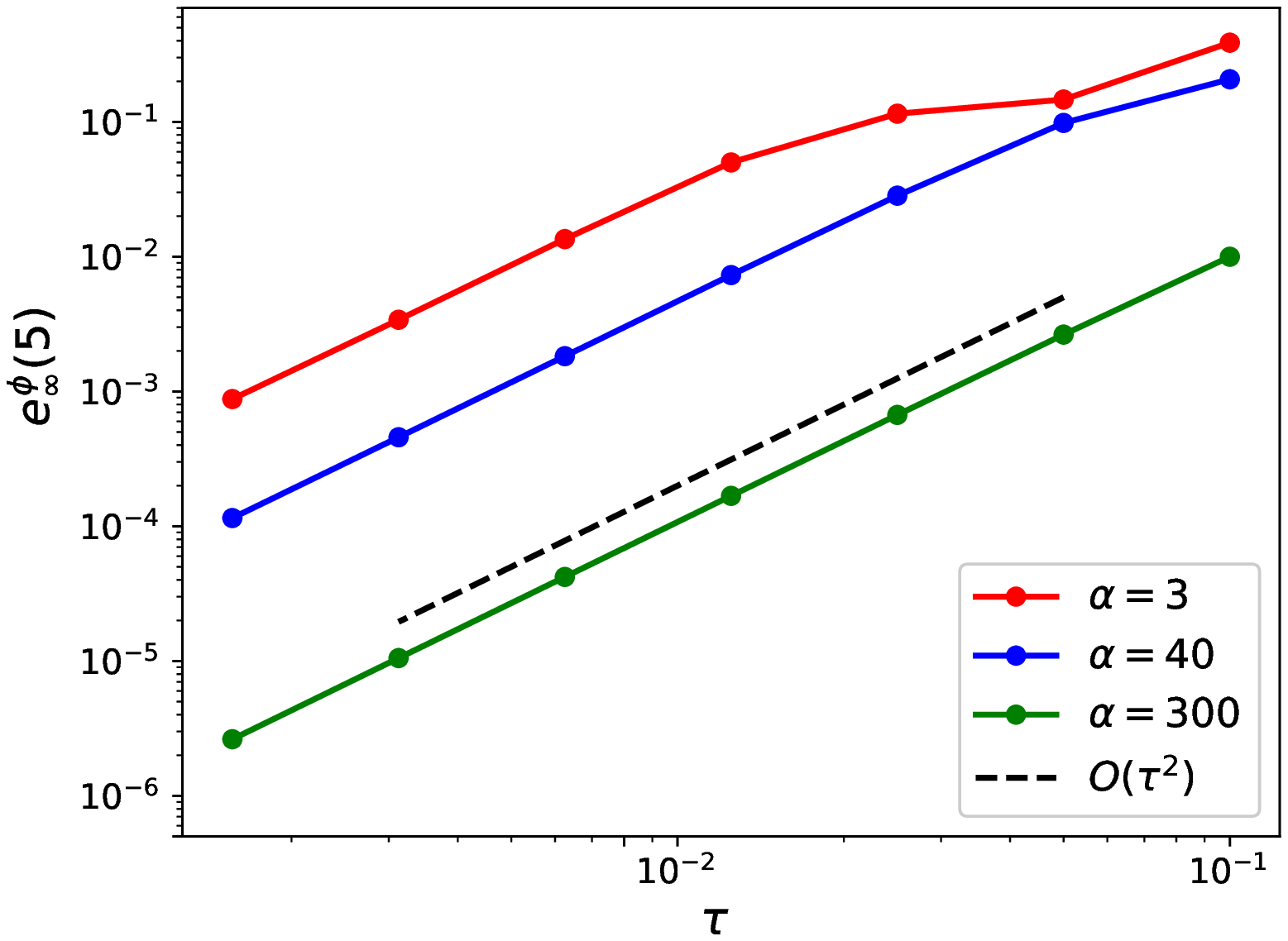}
\end{minipage}}
\hspace{3mm}
\subfigure[Total time $T=n\tau$ required for DSFDN-SI to reach stopping criteria versus $\alpha$ under different time step sizes $\tau$.]{
\begin{minipage}[t]{0.45\linewidth}
\centering
\includegraphics[width=2.9in,height=2.6in]{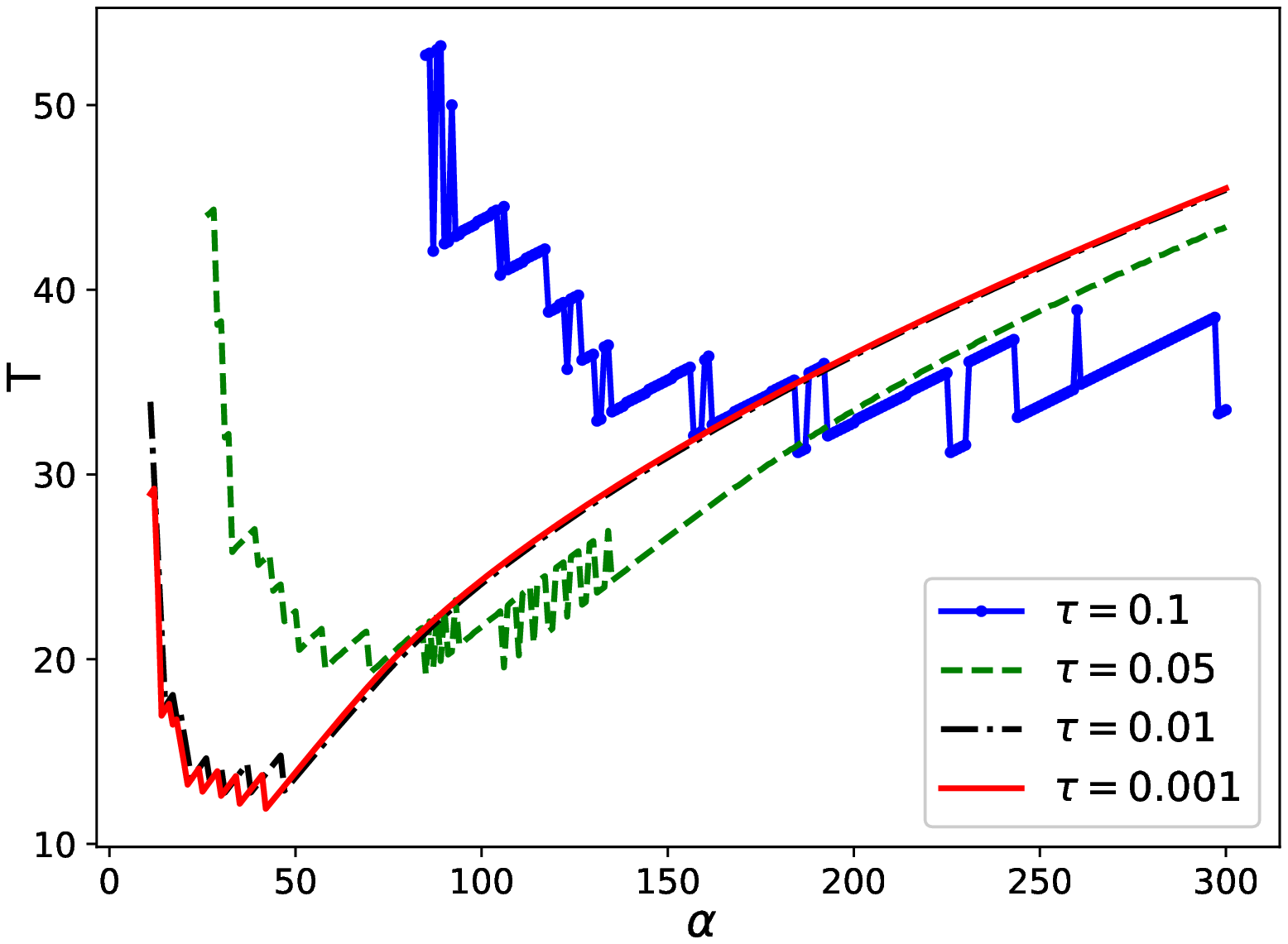}
\end{minipage}}%
\caption{Numerical tests of DSFDN-SI with different $\alpha$ for 1D cases in Example 5.1.}
\label{RespLettFig1}
\end{figure}

{
We first numerically check the impact of damping parameter $\alpha$ on the DSFDN-SI as a showcase example. Fig.~\ref{RespLettFig1}(a) visualizes the temporal error order for approximating $\phi(\cdot,t=5)$ with three different values of $\alpha$. The error function $e_{\infty}^{\phi}\left(t_n\right)=\left\|\phi\left(\cdot, t_n\right)-\phi^{n}\right\|_{\infty}$ at $t=t_n=n\tau$ is introduced here. It appears that numerically larger $\alpha$ gives more stable iterates, and a stable second-order temporal accuracy is achieved for all $\alpha>3$ when $\tau<0.01$. Based on this observation, we fix four different time step sizes, vary $\alpha$ from 3 to 300, and test the total time consumed in order to reach the stopping criteria. Fig.~\ref{RespLettFig1}(b) shows that when the time step size $\tau$ is small (e.g., $\tau=0.01$ or $0.001$), the fastest convergence with respect to $t$ is achieved around $\alpha=40$.  It is believed that this result is close to that of the continuous PDE system of DSFDN since the small time step sizes yield a sufficiently accurate approximation to the original continuous system. On the other hand, we observe in Fig.~\ref{RespLettFig1}(b) that DSFDN-SI achieves the fastest convergence with respect to $t$ when $\alpha$ takes a larger value and a larger step size is used (e.g., $\tau=0.1$). Although the approximation of the continuous flow by DSFDN-SI could be rough for such large time step, as a numerical algorithm, fewer iterations indicate less CPU computational time, which is still appealing as an algorithm for finding the BEC ground states. Therefore, in later examples, we will not exclude to use large time step such as $\tau=0.1$. The experiments here provide intuition for selecting the damping parameter.} 

We then explore the numerical efficiency of DSFDN compared with GFLM. Table~\ref{tab:table1} exhibits the results. For the semi-implicit schemes, when the time step size $0.05\leq\tau\leq1$, the DSFDN-SI only needs approximately one-tenth of the iterations required for the GFLM-BF to reach the stopping condition, and the CPU time consumed is correspondingly shorter. {For the explicit schemes, the maximal size of a stable time step allowed for the GFLM-FE is approximately 0.0006, whereas the DSFDN-FE can choose up to 0.025, which is slightly larger than $\sqrt{0.0006}$. This result corroborates our discussion in Section~\ref{sec:disc-DSFDN} and illustrates that the DSFDN-LF can boost performance by rescaling the time step sizes compared with the gradient flow methods.}

\begin{table}[!t]
\caption{Comparison of numerical results computed by the DSFDN-SI, GFLM-BF, DSFDN-LF and GFLM-FE schemes for the 1D BEC in Example~\ref{5.1}.}
\vspace{2mm}
\label{tab:table1}
\centering
\small\renewcommand{\arraystretch}{1.1}
\begin{tabular}{@{\extracolsep{10pt}}cccccccc@{}}
\hline
Method & $\alpha$ & $\tau$ & Iter & CPU(s) & $E_g$ & $\mu_g$ & $e_{r}^{n}$ \\
\hline
\multirow{4}{*}{DSFDN-SI } & 30000 & 1 & 576 & 0.34 & 15.62475 & 26.01221 & 6.13E-13 \\
& 5000 & 0.5 & 489 & 0.28 & 15.62475 & 26.01221 & 7.26E-13 \\
 & 300 & 0.1 & 564 & 0.31 & 15.62475 & 26.01221 & 5.58E-13 \\
& 100 & 0.05 & 646 & 0.36 & 15.62475 & 26.01221 & 7.77E-13 \\
& 100 & 0.025 & 1196 & 0.63 & 15.62475 & 26.01221 & 9.39E-13 \\
\hline
\multirow{4}{*}{GFLM-BF } & - & 1 & 5804 & 2.25 & 15.62475 & 26.01221 & 9.94E-13 \\
 & - & 0.5 & 5812 & 2.24 & 15.62475 & 26.01221 & 9.95E-13 \\
 & - & 0.1 & 5854 & 2.25 & 15.62475 & 26.01221 & 1.00E-12 \\
 & - & 0.05 & 5912 & 2.29 & 15.62475 & 26.01221 & 9.94E-13 \\
 & - & 0.025 & 6025 & 2.32 & 15.62475 & 26.01221 & 9.95E-13 \\
\hline
\multirow{3}{*}{DSFDN-LF } & 100 & 0.025 & 1137 & 0.46 & 15.62475 & 26.01221 & 4.90E-13 \\
 & 100 & 0.02 & 1422 & 0.62 & 15.62475 & 26.01221 & 4.35E-13 \\
 & 100 & 0.01 & 2844 & 1.13 & 15.62475 & 26.01221 & 8.88E-13 \\
\hline
\multirow{3}{*}{GFLM-FE } & - & 0.0006 & 9317 & 2.14 & 15.62475 & 26.01221 & 9.98E-13 \\
 & - & 0.0004 & 13982 & 3.18 & 15.62475 & 26.01221 & 9.96E-13 \\
 & - & 0.0001 & 55980 & 12.68 & 15.62475 & 26.01221 & 9.99E-13\\
\hline
\end{tabular}
\label{table1}
\end{table}

\begin{figure}[!t]
\hspace*{2em}
\subfigure[Number of iterations to converge for different algorithms with various step sizes $\tau$.]{
\begin{minipage}[t]{0.45\linewidth}
\centering
\includegraphics[width=2.9in,height=2.6in]{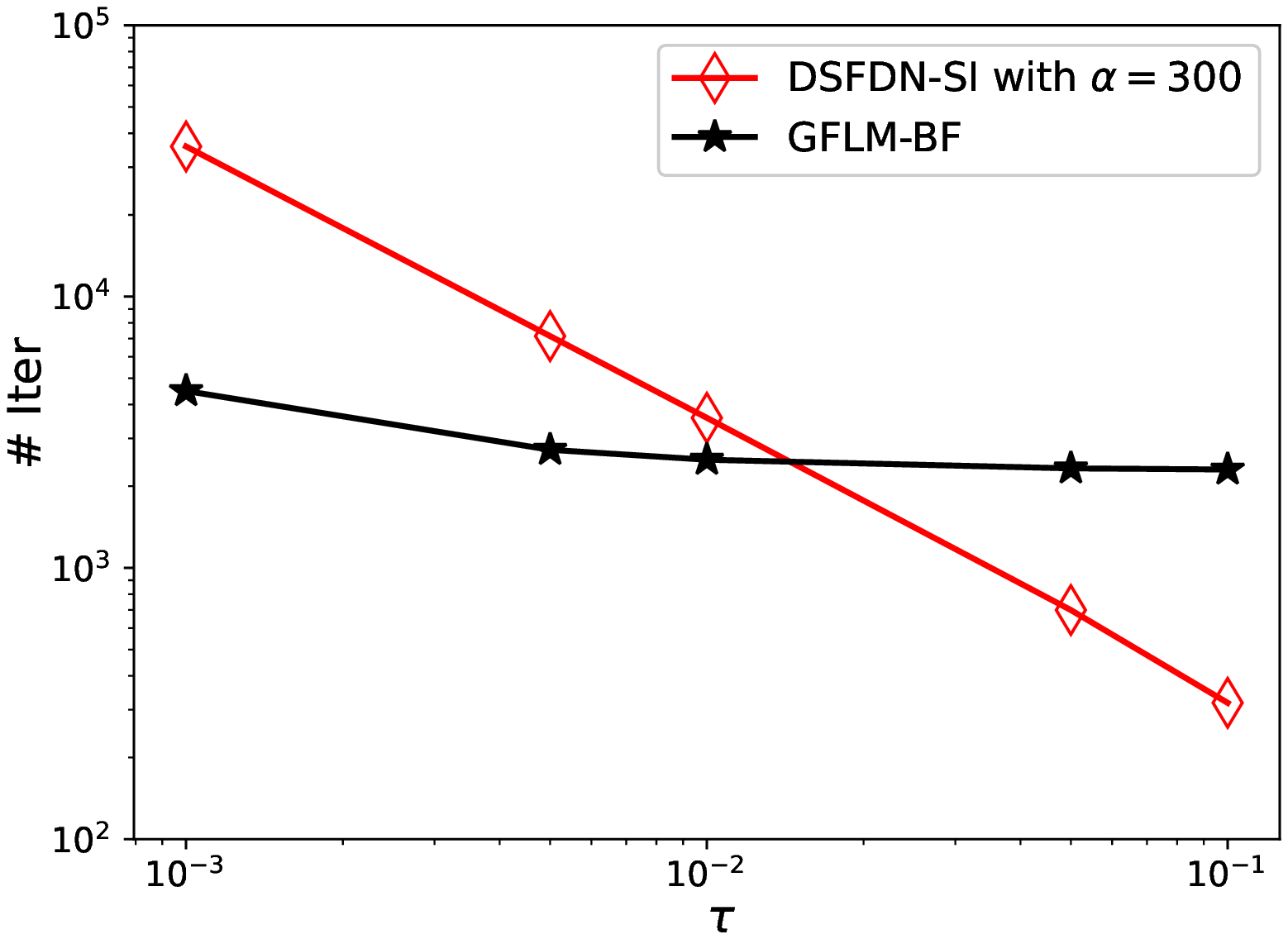}
\label{fig1.a}
\end{minipage}}%
\hspace{3mm}
\subfigure[Errors $e_{\infty}^{\lambda}\left(1\right)$ and $e_{\infty}^{\phi}\left(1\right)$ of the DSFDN-SI scheme and the DSFDN-LF scheme.]{
\begin{minipage}[t]{0.45\linewidth}
\centering
\includegraphics[width=2.9in,height=2.6in]{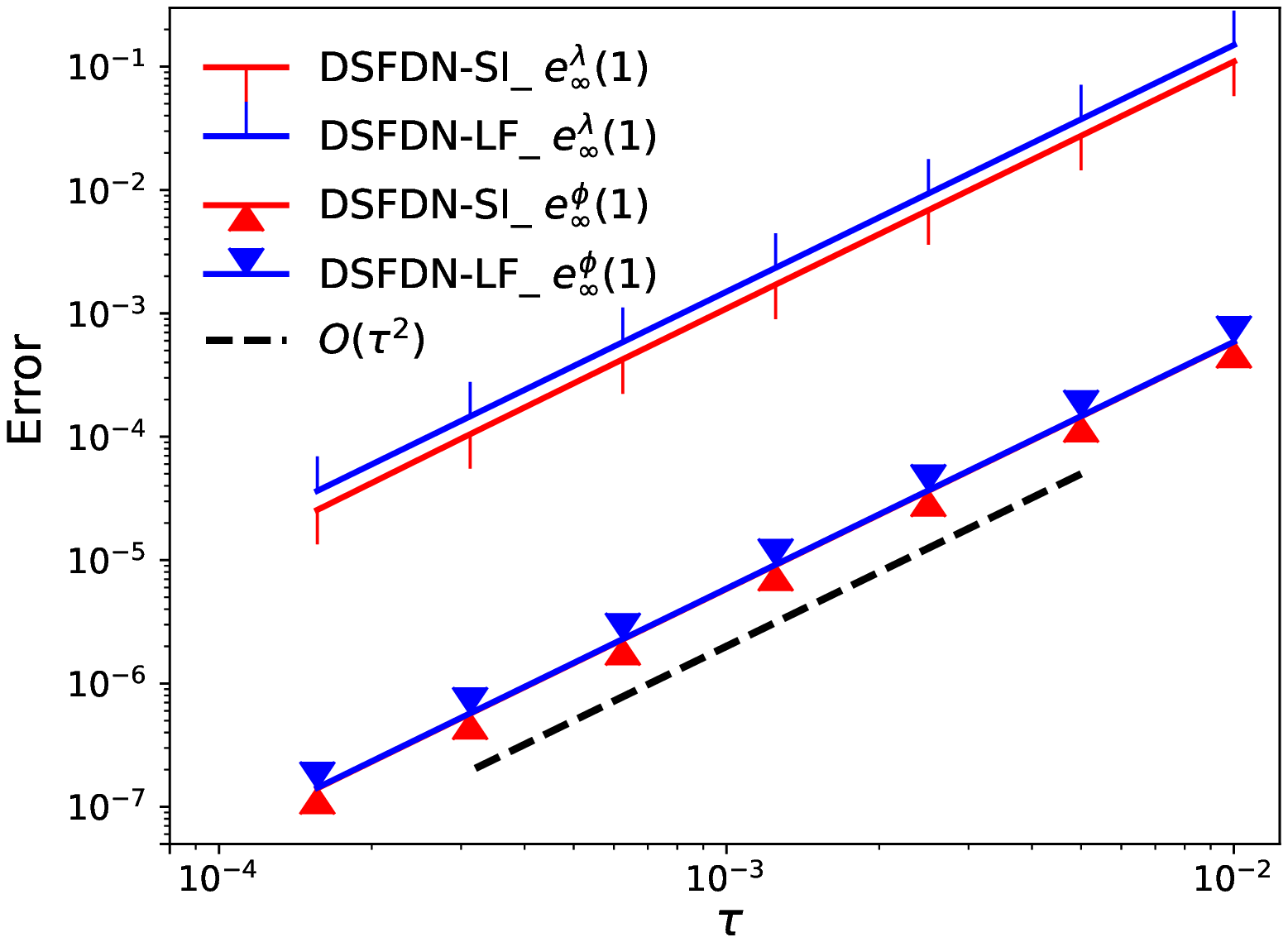}
\label{fig1.b}
\end{minipage}}
\caption{Numerical results of different numerical schemes for 1D cases in Example~\ref{5.1}.}
\end{figure}

Next we verify the temporal discretization schemes described in Section~\ref{sec:disc-DSFDN} holding $O(\tau^2)$ accuracy. It can be seen from Fig.~\ref{fig1.a} that for a fixed $\alpha$, the number of iterations required to converge for the DSFDN-SI is approximately inversely proportional to the time step size $\tau$ when it is suitably small (e.g., $\tau\leq0.1$ for $\alpha=300$ as shown in the figure). This illustrates that the DSFDN-SI, as a second-order accurate temporal discretization scheme, can approach well its continuous form. Fig.~\ref{fig1.b} further visualizes the temporal error order for the DSFDN-SI and DSFDN-LF {by examining the error functions $e_{\infty}^{\lambda}(t_n)=|\lambda_{\phi}\left(t_n\right)-\lambda^n|$ and $e_{\infty}^{\phi}(t_n)=\|\phi\left(\cdot, t_n\right)-\phi^n\|_{\infty}$ at $t=t_n$. Fig.~\ref{fig1.b} shows that both schemes achieve $O(\tau^2)$ time accuracy for approximating $\lambda_{\phi}(t)$ and $\phi(\cdot,t)$ at $t=1$.} It should be noted that the $O(\tau^2)$ time error vanishes at the convergent state, as mentioned in Section~\ref{sec:proposed_flows}-\ref{sec:disc}.

\begin{figure}[!t]
\centering
\hspace*{-3ex}
\includegraphics[width=.9\textwidth,height=.4\textwidth]{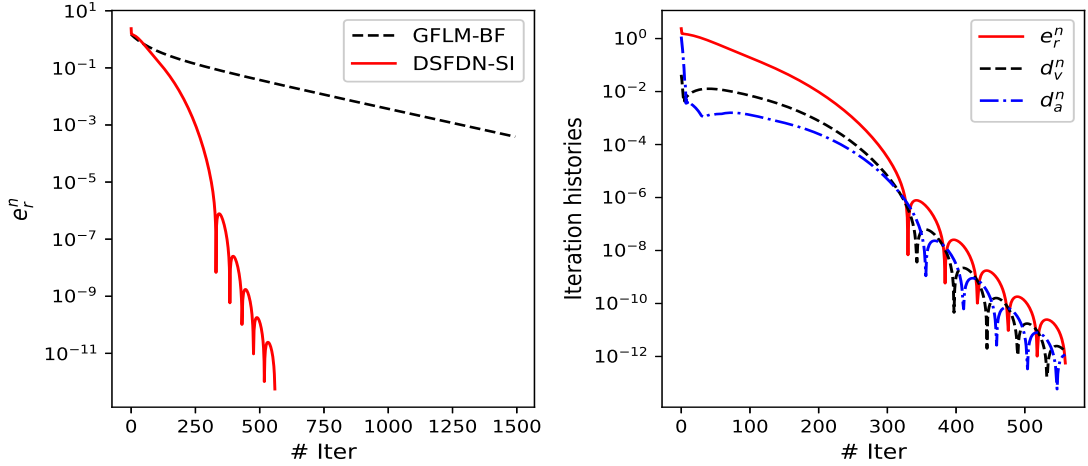}
\caption{Convergence behavior of DSFDN-SI for Example~\ref{5.1}. {\bf Left:} Decaying process comparison for $e^{n}_{r}$ of DSFDN-SI and GFLM-BF. {\bf Right:} Decaying process of various quantities for DSFDN-SI: the maximal residual of Euler-Lagrange equation $e_{r}^{n}$, the discretized velocity $d_{v}^{n}$ and the discretized acceleration $d_{a}^{n}$.}
\label{fig:eg3}
\end{figure}

Lastly, we fix $\tau=0.1$ and $\alpha=300$, and explore the convergence behaviors of the DSFDN-SI scheme by investigating the evolution of various variables. We first compare the decaying phenomenon concerning the maximal residual $e_{r}^{n}$ \eqref{eq:residual-cre}  of DSFDN-SI and GFLM-BF. From left subfigure of Fig.~\ref{fig:eg3}, we observe that for the GFLM-BF scheme, $e_{r}^{n}$ shows a linear exponential downward trend, whereas for the DSFDN-SI scheme, $e_{r}^{n}$ shows a faster oscillating decline. Then we show the change curves of the three quantities $e_{r}^{n}$ \eqref{eq:residual-cre}, $d_{v}^{n}$ \eqref{eq:velocity-cre}, and $d_{a}^{n}$ in DSFDN-SI during the iterative process, where 
\begin{align}
&d_{a}^{n}:=\frac{\big\|\tilde{\phi}^{n+1}-2 \phi^{n}+\phi^{n-1}\big\|_{\infty}}{\tau^2} , \quad n\geq 1.
\end{align} From right subfigure of Fig.~\ref{fig:eg3}, all the three quantities in DSFDN-SI decay with oscillations, and  $e_{r}^{n}$ can be lower than $d_{v}^{n}$ and $d_{a}^{n}$. This oscillation phenomenon is one of the features of inertial dynamics, which is also typically observed in Nesterov's (or Polyak's) accelerated gradient method and the second-order flow systems for convex optimization problems.

In all the subsequent examples in 1D, we always adopt semi-implicit numerical schemes with $\tau=0.1$ and $\alpha=300$.

\begin{figure}[!t]
\centering
\includegraphics[width=.9\textwidth,height=.4\textwidth]{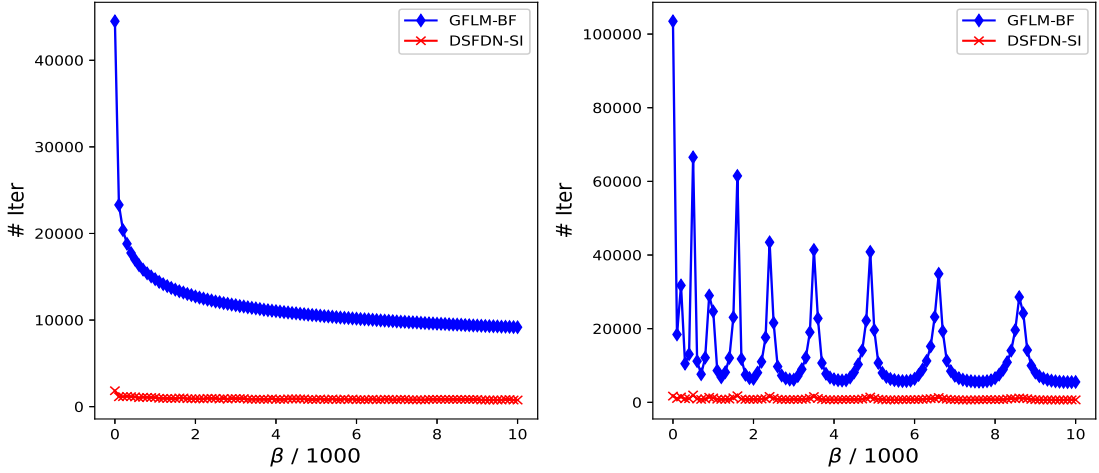}
\caption{Number of iterations of the DSFDN-SI and GFLM-BF schemes to converge for various values of $\beta$ under the harmonic potential (left) and the harmonic-plus-lattice potential (right) in Example~\ref{eg:5.2}.}
\label{fig:eg4}
\end{figure}

\begin{example}\label{eg:5.2}
In this example, we compare the robustness of the DSFDN approach and the GFLM approach in 1D problems. Fig.~\ref{fig:eg4} shows the number of iterations required to converge of the DSFDN-SI and GFLM-BF schemes under different potential functions when $\beta$ varies from $1$ to $10000$. From this figure, we can observe that the DSFDN-SI scheme shows more efficient and robust results than the GFLM-BF scheme for both the harmonic potential \eqref{eq:harmonic-potential} and the harmonic-plus-lattice potential \eqref{eq:harmonic-p-o-l}. In particular, under the harmonic-plus-lattice potential, the different values of $\beta$ have a significant impact on the convergence rate of the GFLM-BF scheme. In contrast, the number of iterations required for the convergence of the DSFDN-SI scheme for different $\beta$ is always in a small range (between 500 and 2000), which indicates an insensitivity of the DSFDN-SI scheme to the nonlinearity strength $\beta$ and shows the strong robustness of our DSFDN approach.
\end{example}

\begin{figure}[!t]
\centering
\includegraphics[width=\textwidth]{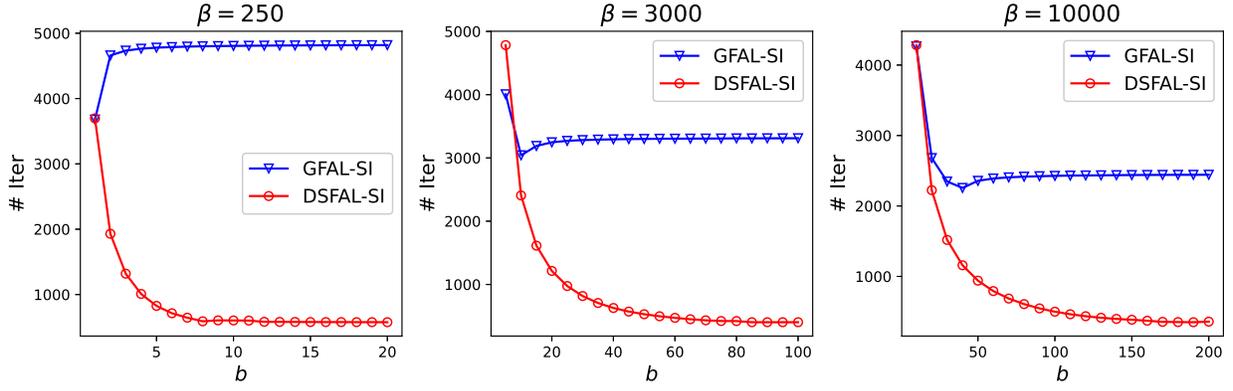}
\caption{Number of iterations of the DSFAL-SI and GFAL-SI schemes to converge for various values of the constant $b$ appeared in the scaling factor $\xi(t)$ in \eqref{xi-form} under different nonlinear strengths $\beta$ in Example~\ref{eg:5.3}.}
\label{fig:eg14}
\end{figure}

\begin{example}\label{eg:5.3}
In this example, we investigate the numerical performance of the DSFAL approach. To this end, we take $V(\mathbf{x})$ as the harmonic potential \eqref{eq:harmonic-potential} and test the convergence results of the DSFAL-SI scheme \eqref{eq:dsfal-si} and the GFAL-SI scheme \eqref{eq:gfal-si} for different $\beta$. Fix $\sigma=100$. Fig.~\ref{fig:eg14} shows the number of iterations to converge against various constant $b>0$ appeared in the scaling factor $\xi(t)$ \eqref{xi-form} when $\beta$ is $250$, $3000$ and $10000$. From Fig.~\ref{fig:eg14} and additional results not shown here for brevity, we find that the constant $b$ affects significantly the overall convergence rate of the DSFAL-SI scheme when it is not too large, and the larger $b$ is, the smaller number of time steps is required to reach the stopping criteria (i.e., the faster the convergence is).
The results in Fig.~\ref{fig:eg14} show that, compared with the GFAL-SI scheme, the DSFAL-SI scheme requires a larger constant $b$ to make the gradient flow of $\chi$ match the convergence speed of the damped second-order flow with respect to $\phi$ and obtain a better overall convergence rate.

On the other hand, Table~\ref{tab:table2} shows a comparison of the DSFAL-SI scheme and the DSFDN-SI scheme when $\beta$ is $250$, $3000$ and $10000$. From Table~\ref{tab:table2}, we observe that the DSFAL-SI scheme with a suitably large $b$ can achieve a comparable numerical performance as the DSFDN-SI scheme in this 1D example.
\end{example}

\begin{table}[!t]
\caption{Numerical comparison of the DSFDN-SI and DSFAL-SI schemes for different $\beta$ in Example~\ref{eg:5.3}.}
\vspace{2mm}
\label{tab:table2}
\centering
\small
\renewcommand{\arraystretch}{1.1}
\begin{tabular}{@{\extracolsep{10pt}}ccccccccc@{}}
\hline
\multicolumn{1}{l}{Method} & $\beta$ & $\tau$ & $b$ & Iter & CPU(s) & $E_g$ & $\mu_g$ & $e_{r}^{n} $\\
\hline
\multirow{3}{*}{DSFDN-SI} & 250 & 0.1 & - & 523 & 0.2898 & 15.62475 & 26.01221 & 9.81E-13 \\
 & 3000 & 0.1 & -  & 434 & 0.2282 & 81.77652 & 136.2867 & 4.71E-13 \\
 & 10000 & 0.1 & - & 379 & 0.2035 & 182.4691 & 304.1114 & 4.69E-13 \\
 \hline
\multirow{3}{*}{DSFAL-SI} & 250 & 0.1 & 20 & 677 & 0.3947 & 15.62475 & 26.01221 & 9.72E-13 \\
 & 3000 & 0.1 & 100 & 436 & 0.2625 & 81.77652 & 136.2867 & 7.35E-13 \\
 & 10000 & 0.1 & 200 & 399 & 0.2343 & 182.4691 & 304.1114 & 9.88E-13\\
 \hline
\end{tabular}%
\end{table}

\subsection{Numerical results and comparisons in 2D examples}
Then we compare different methods in 2D. Unless otherwise specified, the computational domain and the mesh size are chosen to be $\mathcal{D}=[-16,16]^2$ and $h=\frac{1}{32}$, respectively. The following seven types of initial guesses for ground states are used in our numerical computations in rotating cases in 2D \cite{ALT2017JCP,BWM2005CMS}:
\begin{subequations}\label{init}
\begin{align}
\label{init_a} &(a)\ \phi_{a}(\mathbf{x})=\phi_{\mathrm{ho}}(\mathbf{x})=\frac{1}{\sqrt{\pi}} e^{-\left(x^{2}+y^{2}\right) / 2}, && \\ 
\label{init_b} &(b)\ \phi_{b}(\mathbf{x})=\phi_{\text{ho}}^{v}(\mathbf{x})=(x+i y) \phi_{\mathrm{ho}}(\mathbf{x}), &&\hspace{-3em} (\overline{b})\  \phi_{\overline{b}}(\mathbf{x})=\overline{\phi}_{b}(\mathbf{x}),\\
\label{init_c}
 &(c)\ \phi_{c}(\mathbf{x})=\frac{\phi_{\mathrm{ho}}(\mathbf{x})+\phi_{\mathrm{ho}}^{v}(\mathbf{x})}{\left\|\phi_{\mathrm{ho}}+ \phi_{\mathrm{ho}}^{v}\right\|}, &&\hspace{-3em} (\overline{c})\  \phi_{\overline{c}}(\mathbf{x})=\overline{\phi}_{c}(\mathbf{x}),\\
\label{init_d}
 &(d)\ \phi_{d}(\mathbf{x})=\frac{(1-\Omega) \phi_{\mathrm{ho}}(\mathbf{x})+\Omega \phi_{\mathrm{ho}}^{v}(\mathbf{x})}{\left\|(1-\Omega) \phi_{\mathrm{ho}}+\Omega \phi_{\mathrm{ho}}^{v}\right\|}, &&\hspace{-3em} (\overline{d})\ \phi_{\overline{d}}(\mathbf{x})=\overline{\phi}_{d}(\mathbf{x}).
\end{align}
\end{subequations}
The tolerances in stopping conditions \eqref{eq:residual-cre}-\eqref{eq:velocity-cre} and \eqref{eq:stop-pd} are taken as $\varepsilon_v=\varepsilon_r=10^{-10}$ and $\varepsilon_c=10^{-12}$ for all 2D cases.

\begin{example}\label{eg:5.4}
In this example, we illustrate the efficiency of DSFDN-SI through abundant comparison with other methods. Two different problem settings are considered here:
\vspace{-8pt}
\begin{enumerate}[(i)]
\item a fixed $\Omega$ with various $\beta$; 
\vspace{-8pt}
\item a fixed $\beta$ with various $\Omega$.
\end{enumerate}
\end{example}
\vspace{-8pt}

\begin{figure}[!t]
\centering
\includegraphics[width=.9\textwidth,height=.4\textwidth]{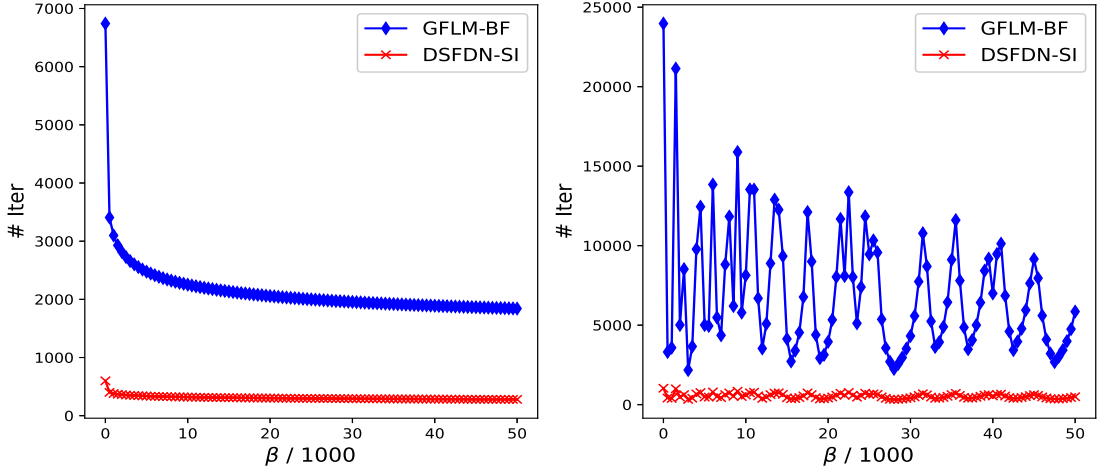}
\caption{Number of iterations of the DSFDN-SI scheme and the GFLM-BF scheme in 2D non-rotating cases for various values of $\beta$ under the harmonic potential (left) and the harmonic-plus-lattice potential (right).}
\label{fig:eg5}
\end{figure}

We first set $\Omega=0$ and vary the value of $\beta$ from 1 to 50000. 
$\phi_0(\mathbf{x})$ given in \eqref{init_tf} is taken as our initial guess for ground state. The left-hand side of Fig.~\ref{fig:eg5} shows the number of iterations to converge against various $\beta$ with $V(\mathbf{x})$ being taken as the harmonic potential \eqref{eq:harmonic-potential} with $\gamma_x=\gamma_y=1$, and the right-hand side shows the situation when $V(\mathbf{x})$ is taken as the harmonic-plus-lattice potential \eqref{eq:harmonic-p-o-l} with $\gamma_x=\gamma_y=1$, $\kappa=25$ and $q_x=q_y=\frac{\pi}{2}$. From Fig.~\ref{fig:eg5}, similar outcomes to 1D cases (see Fig.~\ref{fig:eg3}) are further observed in 2D cases without rotation. The superiority of the DSFDN-SI method in terms of convergence rate compared to the GFLM-BF method is also prominent for 2D cases.

\begin{table}[!t]
\centering
\caption{Numerical results computed by GFLM-BF and DSFDN-SI with different parameters, for 2D rotating BECs in Example~\ref{eg:5.4} with $\beta= 500$ and different $\Omega$. {The maximum iterations allowed for all algorithms are set as 500000.}}
\vspace{2mm}
\label{tab:table3}
\small
\renewcommand{\arraystretch}{1.1}
\begin{tabular}{@{\extracolsep{12pt}}clllllll@{}}
\hline
\multicolumn{1}{l}{$\Omega$} & Method & $\eta$  & $\tau$ & Iter & CPU(s) & Energy & $e_{r}^{n}$ \\
\hline
\multirow{3}{*}{0.5} & {\small GFLM-BF} & - & 0.01 & 500000 & 42839  & 8.032225 & 9.59E-05  \\
 & {\small DSFDN-SI}& $30/t$ & 0.01 & 157880 & 15984 & 8.019671 & 9.21E-11 \\
 & {\small DSFDN-SI}& $50/t$ & 0.1 & 23838 & 2510 & 8.019671 & 9.79E-11 \\
  & {\small DSFDN-SI}& $4$ & 0.1 & 114796 & 11630 & 8.019671 & 9.85E-11 \\
 \hline
\multirow{3}{*}{0.6} & {\small GFLM-BF}& - & 0.01 & 500000 & 42697 & 7.584454 &1.41E-08
  \\
 & {\small DSFDN-SI}& $30/t$ & 0.01 & 66402 & 6222 & 7.584454 & 4.50E-11 \\
 &{\small DSFDN-SI}& $100/t$ & 0.1 & 10290 & 1092 & 7.584454 & 2.29E-11 \\
 \hline
\multirow{3}{*}{0.7} &{\small GFLM-BF}& - & 0.01 & 500000 & 42815 & 6.972631
 & 5.55E-07
 \\
 & {\small DSFDN-SI}& $30/t$ & 0.01 & 86391 & 8960 & 6.972631 & 9.96E-11 \\
 & {\small DSFDN-SI}& $300/t$ & 0.1 & 17211 & 1827 & 6.972631 & 1.40E-11 \\
 \hline
\multirow{3}{*}{0.8} & {\small GFLM-BF} & - & 0.01 & 500000 & 42996 & 6.102827
 & 8.35E-04 \\
 & {\small DSFDN-SI}& $30/t$ & 0.01 & 73738 & 7723 & 6.099744 & 3.60E-11 \\
 & {\small DSFDN-SI}& $200/t$ & 0.1 & 18843 & 2001 & 6.099744 & 7.98E-12 \\
 \hline
\multirow{3}{*}{0.9}& {\small GFLM-BF}& -  & 0.01 & 500000 & 42890 &4.780601  &5.78E-05
  \\
 &{\small DSFDN-SI}& $30/t$ & 0.01 & 254448 & 26307 & 4.777739 & 9.98E-11 \\
 & {\small DSFDN-SI}& $50/t$ & 0.1 & 24773 & 2614 & 4.777739 & 9.99E-11\\
 \hline
\end{tabular}%
\end{table}

Then we consider the rotating BEC in 2D cases with angular velocity $\Omega>0$. Fix $\beta=500$ for different $\Omega$. $V(\mathbf{x})$ is chosen to be the harmonic potential \eqref{eq:harmonic-potential} with $\gamma_x=\gamma_y=1$. The initial guess for ground state is taken as $\phi_{\overline{d}}(\mathbf{x})$ \eqref{init_d}. The maximum number of iterations allowed for all algorithms is set to be $500000$. Our numerical computations show that in this example when $\Omega\geq0.5$, with the same stabilization factor given in \eqref{eq:stabn}, the allowed maximum time step size of the GFLM-BF scheme for different $\Omega$ is approximately 0.01, whereas the DSFDN-SI scheme can take the time step size up to 0.1. Therefore, we only test the situation $\tau =0.01$ for the GFLM-BF scheme, and test two cases of $\tau=0.1$ and $\tau=0.01$ for the DSFDN-SI scheme. We uniformly select $\alpha=30$ in the DSFDN-SI scheme when $\tau=0.01$. For $\tau=0.1$, since the computational cost is less intensive, we test the DSFDN-SI scheme with different $\alpha=$ $50$, $100$, $150$, $200$, $250$, and $300$, and select the optimal result among which.

We make some comments on the work in \cite{OG2020EJP,GOJPA} here. Compared with the DFPM in \cite{OG2020EJP,GOJPA}, the two families of damped second-order flows proposed in this paper adopt completely different ways of handling constraints. Moreover, efficient second-order time-accurate discretization schemes are employed in this paper rather than the symplectic Euler algorithm in \cite{OG2020EJP,GOJPA} with first-order accuracy for the time derivatives. The second-order discretization of the time derivatives can better represent the continuous flows without increasing the computational complexity in comparison with first-order discretizations, thus obtaining faster convergence rates when a suitably larger time step is applied. In addition, numerical experiments in this paper take the general damping coefficient $\eta(t)$ as $\alpha/t$ for some $\alpha\geq3$ instead of simply using constant damping parameters as in \cite{OG2020EJP,GOJPA}. In that regard, we test an example where $\eta$ is constant and list the result in Table~\ref{tab:table3} but using the same discretization scheme as DSFDN-SI. Referring to \cite{OG2020EJP}, we take $\eta=4$, $\tau=0.1$, and test the case of $\Omega=0.5$. Even though the discretization is different from \cite{OG2020EJP,GOJPA}, they approximate the same dynamics when the time step is small enough.

The results in Table~\ref{tab:table3} manifest that the second-order flow method also achieves a compelling advantage over the gradient flow method when computing the rotating BEC ground state. First, GFLM-BF fails to reach the convergence condition in all cases within 500000 iterations, whereas DSFDN-SI with the same time step size $\tau=0.01$ can reach the convergence condition before the threshold, and the numerical ground state energy computed by DSFDN-SI is consistent with the result in \cite{ALT2017JCP, WWB2017JSC}. Moreover, since a larger time step size is available for DSFDN-SI, the advantage of its convergence rate is more promising when $\tau=0.1$. Lastly, for the case where $\eta$ is constant, the result of its convergence rate is poor than that of $\eta=\alpha/t$.


\begin{example}\label{eg:5.5}
In this example, we explore the effect of DSFAL-SI computing rotating BEC ground state with setting (ii) in Example~\ref{eg:5.4}. Different from the (non-rotating) case in 1D in Example~\ref{eg:5.3}, the rotational term, especially with high rotating speed, makes the problem more difficult to solve in this case. A larger scaling factor usually results in less robustness of the computational algorithm. We take the constant $b$ in \eqref{xi-form} as $b=0.8$ for $\Omega=0.9$ and $b=1$ for the rest cases in this example.
\end{example}
\vspace{-8pt}

\begin{table}[!t]
\centering
\caption{Numerical results of 2D rotating BECs computed by DSFAL-SI with different augmentation parameter $\sigma$ in Example~\ref{eg:5.5} for $\beta= 500$ and different $\Omega$.}
\vspace{2mm}
\label{tab:table4}
\small
\renewcommand{\arraystretch}{1.1}
\begin{tabular}{@{\extracolsep{15pt}}llllllll@{}}
\hline
$\Omega$ 
& $\sigma$ & $\tau$  & Iter & CPU(s) & Energy & $e_{r}^{n}$ & $e_{c}^{n}$\\
\hline
0.5 
& 200 & 0.1 & 10786 & 1152 & 8.019671 & 1.08E-11 &1.20E-14\\
0.6 
& 100 & 0.1 & 11055 & 1178 & 7.584454 & 9.98E-12 &2.79E-14\\
0.7 
& 100 & 0.1 & 13738 & 1366 & 6.972631 & 9.60E-12 &3.38E-14\\

0.8 
& 100 & 0.1 & 20048 & 2102 &6.099744  & 6.47E-12 &9.33E-14\\

0.9 
&50 & 0.1 & 47958 & 4756 & 4.777739 & 3.57E-11 &5.67E-13\\
\hline
\end{tabular}%
\end{table}

\begin{figure}[!t]
\centering
\subfigure[Time evolution of the density $|\phi|^2$ by DSFAL-SI (top row) and DSFDN-SI (bottom row).]{
\begin{minipage}[c]{0.95\textwidth}
\centering
\includegraphics[width=.95\textwidth]{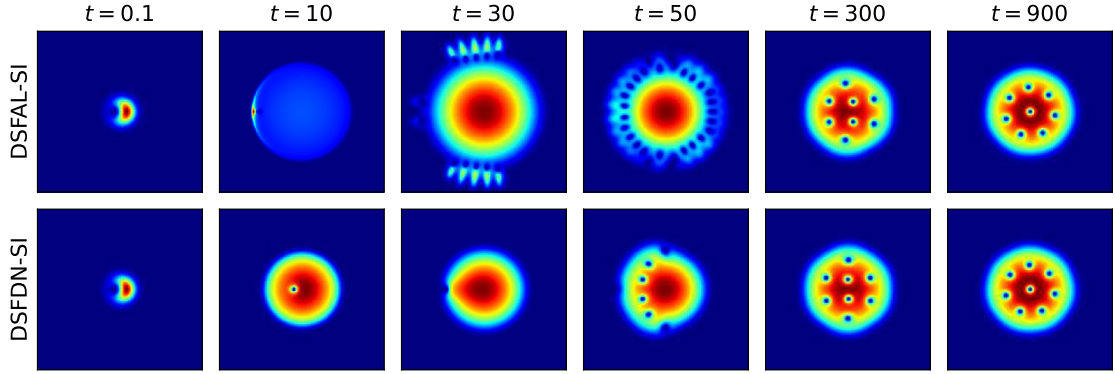}
\end{minipage}
\label{fig:eg7a}
}
\subfigure[GP energy changing during the time evolution of DSFAL-SI and DSFDN-SI. In axis labels, $E(t)$ represents the value of GP energy at time t, and $E_g$ represents the value of the lowest energy during the evolution of DSFDN-SI.]{
\begin{minipage}[c]{0.95\textwidth}
\centering
\includegraphics[width=\textwidth,height=.4\textwidth]{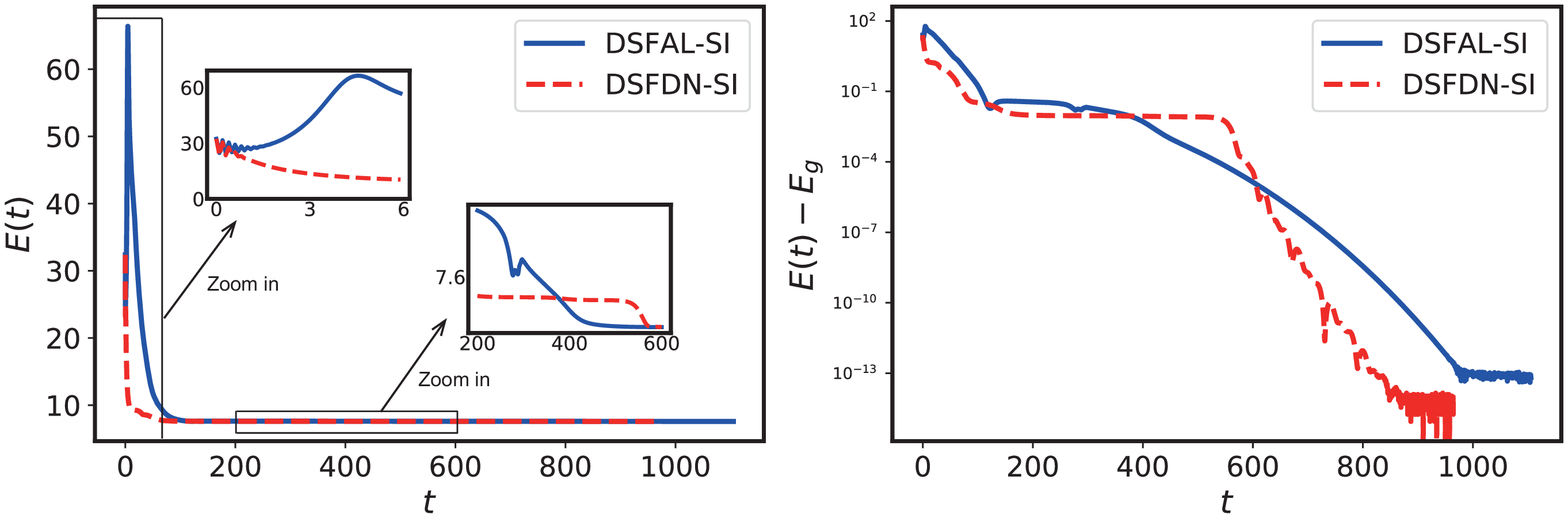}
\end{minipage}
\label{fig:eg7b}
}
\caption{Evolution of DSFAL-SI and DSFDN-SI in 
Example~\ref{eg:5.5} when $\Omega=0.6$.}
\end{figure}

To compare with DSFDN-SI, we use the same initial value and time step for DSFAL-SI. Here we uniformly take $\alpha=200$, and test different augmentation coefficient $\sigma=$ $50$, $100$, and $200$. The optimal results among which is shown in Table~\ref{tab:table4}. Note that the converged energies computed by DSFAL-SI are the same as that of DSFDN-SI, with a comparable convergence rate. However, the evolution process of the two methods from the initial value to the converged solution is significantly different. Fig.~\ref{fig:eg7a} depicts the fictitious time evolution of the density $|\phi|^2$ computed by these two methods for the case $\Omega=0.6$ (Note that we intercepted the figures obtained from the computational domain, to exhibit the evolution process more clearly). It can be seen that throughout the evolution process computed by DSFDN-SI, the contour map of density $|\phi|^2$ always presents a relatively regular image with a relatively fixed size of the bright color (value is larger than 0) area. However, for DSFAL-SI that does not need to satisfy the constraint during the evolution, the contour maps of $|\phi|^2$ have very irregular shapes in the early stages and tend to be regular gradually after a certain time. Eventually, both methods converge to the same state. Fig.~\ref{fig:eg7b} shows the energy changing during the evolution. The energy computed by DSFDN-SI decreases steadily as the fictitious time progresses, while the energy computed by DSFAL-SI increases rapidly at the beginning and then immediately decreases to the same level as DSFDN-SI (see Table~\ref{tab:table3}).

\subsection{{2D problems with fast rotation}}
In this subsection, we apply the DSFDN-SI scheme to rotating BECs with strong interaction and high-speed rotation, which are more related to practical physical problems. A difficulty of this kind of problems is that it is not only highly nonlinear but also highly non-convex, which usually requires expensive computational costs. {We compare the ground state energies computed by DSFDN-SI with those computed by the regularized Newton method (RNM) in \cite{WWB2017JSC}. In addition, the methods based on gradient flow are not considered here because they are too costly for these examples, which may take several days.}

\begin{figure}[!h]
\centering
\vspace{2ex}
\includegraphics[width=1\textwidth]{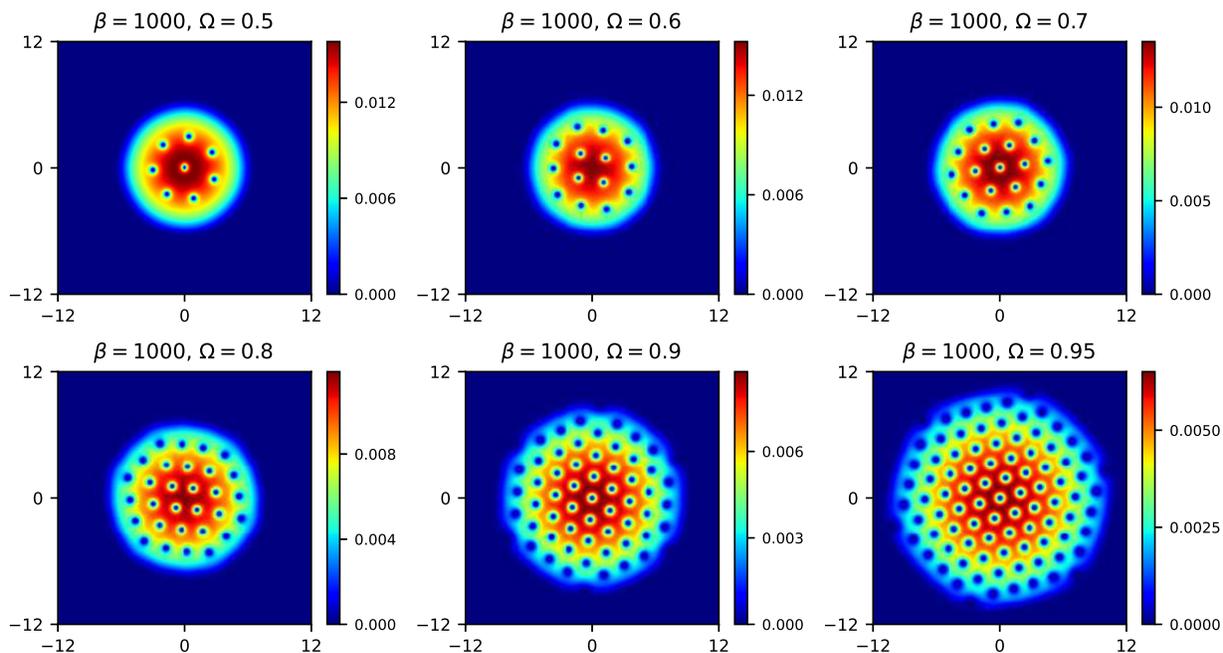}
\caption{Contour plots of $|\phi_g (\mathbf{x})|^2$ corresponding to the lowest energy levels (with  “$\star$” sign) listed in Table~\ref{tab:table5}.}
\label{fig:eg8}
\end{figure}

\begin{example}\label{eg:5.6}
In this example,  we still take $V(\mathbf{x})$ as the harmonic potential \eqref{eq:harmonic-potential} as Example~\ref{eg:5.4} but with a larger value of $\beta$, which increases the nonlinearity of the problem. Here we fix $\beta=1000$, set the computational domain as $[-12,12]^2$, and vary $\Omega$ from 0.5 to 0.95. It is worth noting that for the rotating BEC with harmonic potential, the ground state exists when the angular velocity satisfies $|\Omega|<\min\{\gamma_x,\gamma_y\}=1$ \cite{BWM2005CMS}, thus $\Omega=0.95$ is a considerably high angular rotation speed for this case. The result from \cite{WWB2017JSC} is invoked for comparison, and the cascadic multigrid method (see Remark \ref{multigrid}) used in \cite{WWB2017JSC} is adopted in this example, beginning with the coarsest grid $(2^4+1)\times(2^4+1)$ and ending with the finest grid $ (2^8+1)\times(2^8+1)$. The tolerance parameters are set as $\varepsilon_v=\varepsilon_r=10^{-10}$ for the finest grid and $\varepsilon_v=\varepsilon_r=10^{-6}$ for the rest coarser grids. Table~\ref{tab:table5} lists the converged energy computed by DSFDN-SI with the corresponding CPU time consumed and the parameter selection for various initial data in \eqref{init} and different $\Omega$. The one with the lowest energy among results with different initial values for each given $\Omega$ is marked by a “$\star$” sign. Fig.~\ref{fig:eg8} shows the contour plots $|\phi(\mathbf{x})|^2$ of the convergent solution with lowest energy. It is observed from Fig.~\ref{fig:eg8} that the vortices exhibit a lattice structure at high rotational speeds, forming the so-called Abrikosov lattice (see, e.g., \cite{BC2013KRM,DP2017SISC,Fetter2009RMP}). 

Compared with result provided in \cite[Table~7]{WWB2017JSC}, our results show that when $\Omega=0.5$, $0.6$, $0.7$ and $0.8$, the lowest converged energies computed by DSFDN-SI are consistent with that obtained by RNM \cite{WWB2017JSC}, and when $\Omega=0.9$ and $0.95$, we can get a convergent solution with lower energy via DSFDN-SI.
\end{example}

\begin{table}[!h]
\centering
\caption{Converged GP energy obtained numerically by DSFDN-SI combined with a multigrid technique using different initial values for $\beta=1000$ and different $\Omega$ in Example~\ref{eg:5.6} (with lowest energy for each $\Omega$ marked by ``$\star$").}
\vspace{2mm}
\label{tab:table5}
\small
\renewcommand{\arraystretch}{1.1}
\begin{tabular}{@{\extracolsep{12pt}}lllllll@{}}
\hline
$\Omega$ & 0.5     & 0.6     & 0.7    & 0.8    & 0.9    & 0.95   \\
   \hline
$(a)$ & $11.0954^{\star}$ & 10.4464 & 9.5289 & 8.2610 & 6.3608 & $4.8822^{\star}$ \\
$(b)$ & 11.1369 & 10.4392 & 9.5289 & 8.2633 & 6.3617 & 4.8822 \\
$(\overline{b})$ & 11.1054 & 10.4392 & 9.5289 & 8.2628 & 6.3607 & 4.8822 \\
$(c)$ & 11.1054 & 10.4392 & 9.5301 & 8.2610 & 6.3607 & 4.8822 \\
$(\overline{c})$ & 11.1054 & 10.4392 & 9.5283 & 8.2637 & 6.3607 & 4.8822 \\
$(d)$ & 11.1054 & 10.4392 & 9.5301 & 8.2637 & 6.3607 & 4.8822 \\
$(\overline{d})$ & 11.1054 & $10.4392^{\star}$ & $9.5283^{\star}$ & $8.2610^{\star}$ & $6.3601^{\star}$ & 4.8822\\
\hline
\end{tabular}
\end{table}

\begin{figure}[!t]
\centering
\vspace{2ex}
\includegraphics[width=1\textwidth]{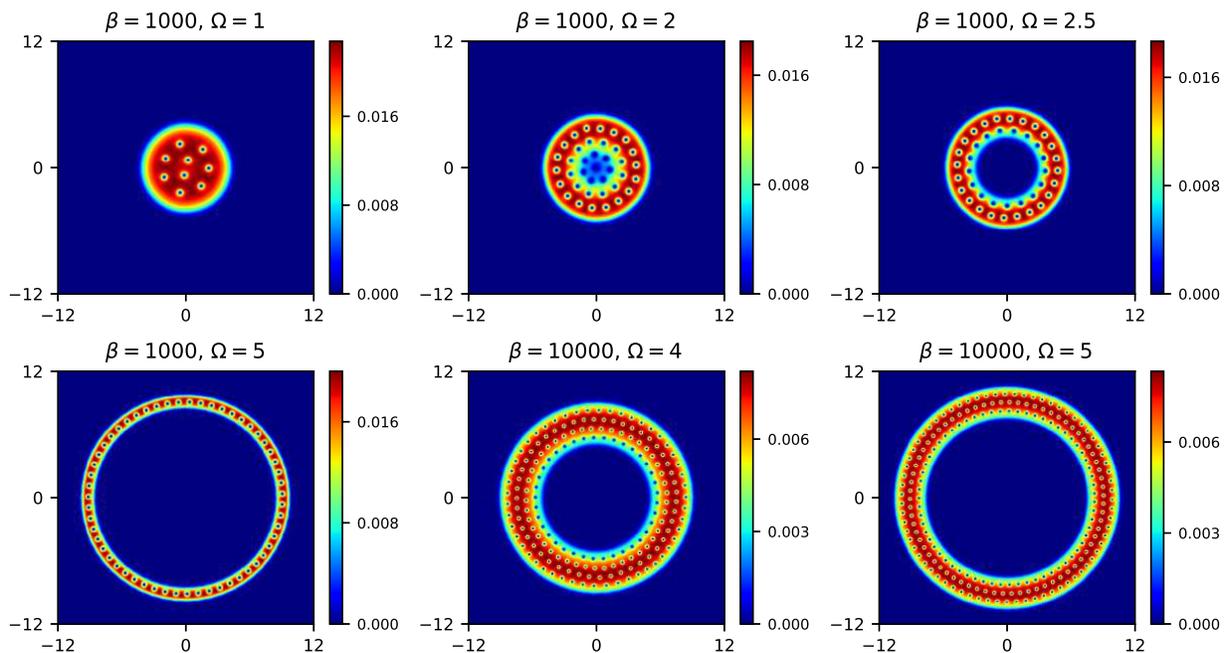}
\caption{Corresponding contour plots of the density function $|\phi_g (\mathbf{x})|^2$ of {Table~\ref{tab:table6}}}
\label{fig:eg9}
\end{figure}

\begin{example}\label{eg:5.7}
In this example, we take $V(\mathbf{x})$ as the harmonic-plus-quartic potential \eqref{eq:harmonic-p-q} with $ \gamma_x= \gamma_y = 1$, {$\theta = 1.2$} and $\kappa = 0.3$. The computational domain is set as $\mathcal{D}=[-12,12]^2$, and the tolerance parameters are chosen as $\varepsilon_r=\varepsilon_v=10^{-7}$. We uniformly take $\phi_{a}(\mathbf{x})$ \eqref{init_a} as the initial value in this example. First, we consider cases tested in \cite{WWB2017JSC} that $\beta=1000$ with different $\Omega=1.0$, $2.0$, $2.5$ and $5.0$, then we compute cases $\beta=10000$ with $\Omega=4.0$ and $5.0$. Spatial mesh size $h$, time step $\tau$, damping parameter $\alpha$, convergent energy, the number of iterations and CPU time consumed are listed in Table~\ref{tab:table6}. Moreover, Fig.~\ref{fig:eg9} shows the contour plot of the density function $|\phi|^2$ of convergent solutions for all these six cases.

From Fig.~\ref{fig:eg9}, it can be observed that the computed solutions are all vortex-symmetric, and a giant hole appeared in the middle for high rotational speed cases. This phenomenon is in line with the characteristics of the ground sate of fast rotating BECs with harmonic-plus-quartic potential as observed in existing experiments \cite{Cornell2003giant} and numerical simulations \cite{aftalion2004giant}, and the spatial division applied here is sufficient to visualize the vortices. Meanwhile, the CPU time consumed in these six cases is acceptable, whereas there were few reports of such problems using gradient flow methods due to the high computational cost. Additionally, comparing our results for $\beta=1000$ with the results in \cite[Table~9]{WWB2017JSC}, it shows that for the case of $\Omega=1.0$, the convergent energy computed by DSFDN-SI is consistent with the one obtained by RNM, and for $\Omega=2.0, 2.5$ and $5.0$, lower energy solutions can be obtained using DSFDN-SI compared to RNM. This result and the one in Example~\ref{eg:5.6} both demonstrate the ability of the second-order flow methods for computing lower-energy solutions of fast rotating BECs. 
\end{example}

\begin{table}[!t]
\centering
\caption{Numerical results computed by DSFDN-SI with different parameters for strong interacting and fast rotating BECs in Example~\ref{eg:5.7}. }
\vspace{2mm}
\label{tab:table6}
\small
\renewcommand{\arraystretch}{1.1}
\begin{tabular}{@{\extracolsep{12pt}}lcccccc@{}}
   \specialrule{0.05em}{0pt}{2pt}
$\beta$ & \multicolumn{4}{c}{1000} & \multicolumn{2}{c}{10000} \\
\cline{2-5} \cline{6-7}
$\Omega$ & 1 & 2 & 2.5 & 5 & 4 & 5  \\
   \specialrule{0.05em}{2pt}{2pt}
$h$ & $24/2^8$ & $24/2^8$ & $24/2^8$ & $24/2^9$ & $24/2^9$ & $24/2^9$ \\
$\tau$ & 0.1 & 0.1 & 0.03 & 0.03 & 0.03 & 0.03 \\
$\alpha$ & 100 & 100 & 40 & 40 & 40 & 40 \\  
Iter & 294749 & 141487 & 48764 & 30037 & 101535 & 289436 \\
CPU(s) & 6522 & 3207 & 1115 & 3211 & 11472 & 34506 \\
Energy & 12.4820 & $-2.3432$ & $-21.7770$ & $-513.7313$ & $-167.2832$ & $-476.5897$ \\
\specialrule{0.05em}{2pt}{0pt}
\end{tabular}
\end{table}

\

\section{Concluding remarks}
\label{sec:concl}
This paper proposed two second-order damped hyperbolic flows to approach the ground states of rotating Bose-Einstein condensates (BECs), namely the damped second-order flow with discrete normalization (DSFDN) and the damped second-order flow based on augmented Lagrangian (DSFAL). The DSFDN extends the framework of the popular normalized gradient flow methods in the literature to second-order flows and redesigns the Lagrange multiplier to always keep track the normalization constraint in the continuous level. The DSFAL is generated from the augmented Langrangian framework that does not need to satisfy the constraint at every moment but only at the end of convergence. Both methods have their advantages. The DSFDN method is easier to implement for the problems with one normalization constraint, such as the single-component BEC under the normalization or mass constraint focused in this paper. However, it is not easy to directly generalize the idea to general multi-component BEC systems with multiple constraints, e.g., the high-spin BEC under the two constraints on the total mass and magnetization. In that regard, DSFAL seems to be more flexible to problems with multiple constraints, though it has more parameters to tune. Nevertheless, both proposed flows lead to novel and interesting topics on second-order hyperbolic type partial differential equations (PDEs) or systems. Theoretical and numerical investigation of these PDEs would give us a deeper understanding, while we have concentrated on the numerical simulations and their applications in computing the ground states of rotating BECs via numerically discretizing the PDE systems. Detailed discretization methods are presented, which generate a series of new algorithms. Extensive one- and two-dimensional numerical examples show that the proposed second-order flow methods are capable of accurately and efficiently computing the ground states of non-rotating and rotating BECs. In particular, the second-order flow methods can achieve better performance in terms of convergence rate and robustness in comparison with their first-order counterparts, i.e., normalized gradient flows. It is not surprising that the DSFDN and the DSFAL often have quite different behaviors during the convergence process, even though they eventually converge to the same state in our experiments. We have observed that numerically the DSFDN-SI, a stabilized semi-implicit discretization of the DSFDN, is an efficient algorithm for dealing with strongly nonlinear and non-convex problems. Comparing the results of existing numerical algorithms based on Newton's method in the literature manifests that the second-order flow method can compute solutions with lower energy.

We believe that the current paper initializes an interesting research topic for further investigation, as the study here is mostly algorithmically and experimentally oriented, and many theoretical questions are open. It would be worth having analysts' efforts on the well-posedness, regularity of solutions of the proposed second-order flows, and also the asymptotic behavior of their trajectories. On the other hand, the convergence and convergence rates of second-order flows for such constrained non-convex minimization problems are also of high interests, both theoretically and numerically.

We also find that our current work might have connections to the so-called phenomenological damping (PD) \cite{ChoMorBur98} of which the original idea is traced back to Pitaevskii in 1958. The mechanism of damping appeared to be the common characteristic of both methods. The difference is that our work is based on some artificial mechanical dynamics from optimization point of view  while PD was directly rooted in the GPE based on quantum mechanics. It is not clear how PD performs numerically in comparison with gradient flow type methods or our second-order flow approaches. It would be interesting to have some deeper understanding on these two types of damping dynamics reflecting different perspectives.

\section*{Acknowledgements}
\addcontentsline{toc}{section}{Acknowledgments}
The work of H. Chen and G. Dong was supported by NSFC grant 12001194.
The work of W. Liu was supported by NSFC grant 12101252, the International Postdoctoral Exchange Fellowship Program No. PC2021024 and the Guangdong Basic and Applied Basic Research Foundation grant 2022A1515010351.
The work of H. Chen and Z. Xie was supported by NSFC grant 12171148.
Part of the work was done when G. Dong was employed by Humboldt University of Berlin and visiting the Mathematical Research Center at Hunan Normal University and when W. Liu was visiting the Department of Mathematics at National University of Singapore.

\appendix

\section{Normalized gradient flow approaches}\label{sec:first-order}
In this section, we describe normalized gradient flow methods for computing the ground state of rotating BECs defined in \eqref{eq:gsdef}. In particular, we extend the GFLM approach proposed in \cite{LC2021SISC,CL2021JCP} to the rotating BEC model and present its typical temporal discretizations.

\subsection{Brief review of conventional normalized gradient flows in rotational frame}
The basic idea of the GFDN is to apply the gradient flow to an unconstrained minimization for the GP energy functional $E(\phi)$ \eqref{eq:GPErot-energy} and then pull the solution back to the normalization constraint manifold $S$. More precisely, setting $t_n=n\tau$ for $n=0,1,\ldots$ with $\tau>0$ a time step length, the GFDN for a rotating BEC reads \cite{BC2013KRM,BD2004SISC,BWM2005CMS,ZZ2009CPC}
\begin{equation}\label{eq:gfdn}\left\{
\begin{aligned}
& \dot{\phi}=-\frac{\delta E(\phi)}{\delta\overline{\phi}}=\frac{1}{2} \Delta\phi-V \phi-\beta|\phi|^{2} \phi+\Omega L_z\phi, \quad\mathbf{x}\in\mathbb{R}^d, \; t\in[t_{n},t_{n+1}), \\
& \phi(\mathbf{x}, t_{n+1}) := \phi(\mathbf{x}, t_{n+1}^{+})=\frac{\phi(\mathbf{x}, t_{n+1}^{-})}{\|\phi(\cdot, t_{n+1}^{-})\|}, \quad\mathbf{x}\in\mathbb{R}^d, \; n\geq0, \\
& \phi(\mathbf{x}, 0)=\phi_{0}(\mathbf{x}), \quad\mathbf{x}\in\mathbb{R}^d,
\end{aligned}\right.
\end{equation}
where $\phi(\mathbf{x}, t_n^{\pm})=\lim_{t\to t_n^{\pm}} \phi(\mathbf{x},t)$ and $\phi_0$ is an initial guess for the ground state satisfying the normalization condition, i.e., $\|\phi_0\|=1$. 
In fact, the GFDN \eqref{eq:gfdn} can be viewed as a first-order splitting for the following CNGF \cite{BD2004SISC,BWM2005CMS}:
\begin{equation}\label{eq:cngf}\left\{
\begin{aligned}
& \dot{\phi}=\frac{1}{2} \Delta\phi-V \phi-\beta|\phi|^{2} \phi+\Omega L_z\phi+\mu_{\phi}(t)\phi, \quad \mathbf{x}\in\mathbb{R}^d,\; t\geq0, \\
& \phi(\mathbf{x}, 0)=\phi_{0}(\mathbf{x}), \quad \mathbf{x}\in\mathbb{R}^d,
\end{aligned}\right.
\end{equation}
with $\|\phi_0\|=1$ and $\mu_{\phi}(t)=\frac{1}{\|\phi(\cdot,t)\|^2}\int_{\mathbb{R}^d}\left(\frac12|\nabla\phi|^2+V|\phi|^2+\beta|\phi|^4-\Omega\overline{\phi}L_z\phi\right)\mathrm{d}\mathbf{x}$. It is proved that the CNGF \eqref{eq:cngf} conserves normalization-constraint and is energy-diminishing \cite{BWM2005CMS}, i.e.,
\begin{align}
\|\phi(\cdot,t)\|^2\equiv\|\phi_0\|^2=1,\quad
\frac{\mathrm{d}}{\mathrm{d}t}E(\phi(\cdot,t))=-2\|\dot{\phi}(\cdot,t)\|^2,\quad \forall t\geq0.
\end{align}
Formally, as $t\to+\infty$, the CNGF \eqref{eq:cngf} will converge to the ground state provided that the initial data $\phi_0$ is properly selected \cite{BWM2005CMS}.

By a similar discussion in \cite{LC2021SISC} for non-rotating BEC case, one can see that the GFDN \eqref{eq:gfdn} suffers from the inaccuracy for finite time step $\tau>0$ due to the $O(\tau)$ operator-splitting error. Although the temporal error of the GFDN \eqref{eq:gfdn} could be eliminated by some special discretizations such as the linearized backward Euler scheme in single-component case \cite{LC2021SISC}, the limitation of the further discretization for the GFDN is inconvenient in practice. We refer to \cite{LC2021SISC} for more details.

\subsection{Normalized gradient flow with Lagrange multiplier for rotating BECs}
Following the work in \cite{LC2021SISC,CL2021JCP}, we propose the GFLM as a modified GFDN to compute the ground state of a rotating BEC as
\begin{equation}\label{eq:gflm}\left\{
\begin{aligned}
&\dot{\phi} = \frac{1}{2} \Delta\phi-V \phi-\beta|\phi|^2 \phi+\Omega L_z\phi + \mu_\phi(t_n)\phi, \quad\mathbf{x}\in\mathbb{R}^d, \; t\in[t_n,t_{n+1}), \\
& \phi(\mathbf{x},t_{n+1}) := \phi(\mathbf{x},t^+_{n+1}) = \frac{\phi(\mathbf{x},t^-_{n+1})}{\|\phi(\cdot,t^-_{n+1})\|},\quad \mathbf{x}\in\mathbb{R}^d,\; n\geq0,  \\
& \phi(\mathbf{x}, 0)=\phi_{0}(\mathbf{x}), \quad \mathbf{x}\in\mathbb{R}^d,
\end{aligned}\right.
\end{equation}
with $\|\phi_0\|=1$ and $\mu_{\phi}(t_n)=\mu(\phi(\cdot,t_n))=E(\phi(\cdot,t_n))+\frac{\beta}{2}\int_{\mathbb{R}^d} |\phi(\mathbf{x},t_n)|^4\mathrm{d}\mathbf{x}$. Apparently, the only difference between the GFLM \eqref{eq:gflm} and the GFDN \eqref{eq:gfdn} is that the former contains an additional Lagrange multiplier term $\mu_\phi(t_n)\phi$ in the gradient flow part. On the other hand, the GFLM \eqref{eq:gflm} can also be viewed as an approximation of the CNGF \eqref{eq:cngf}.

Thanks to the introduction of the Lagrange multiplier term in \eqref{eq:gflm}, the temporal discretization for the GFLM \eqref{eq:gflm} is very flexible. Setting $\phi^0(\mathbf{x})=\phi_0(\mathbf{x})$, we present the following two simple temporal discretization schemes for the GFLM \eqref{eq:gflm}:
\begin{itemize}
\item the forward Euler scheme (GFLM-FE),
\begin{align}\label{eq:gflm-fe}
 \frac{\tilde{\phi}^{n+1}-\phi^n}{\tau} = \left(\frac{1}{2} \Delta -V -\beta|\phi^n|^2 + \Omega L_z + \mu^n\right)\phi^n,
\end{align}
\item the backward-forward Euler scheme (GFLM-BF),
\begin{align}\label{eq:gflm-bf}
\frac{\tilde{\phi}^{n+1}-\phi^n}{\tau} = \left(\frac{1}{2} \Delta -\vartheta^n\right)\tilde{\phi}^{n+1} + \Big(\vartheta^n-V-\beta|\phi^n|^{2}+ \Omega L_z + \mu^n\Big) \phi^n,
\end{align}
\end{itemize}
both followed by a normalization step as $\phi^{n+1} = \tilde{\phi}^{n+1}/\|\tilde{\phi}^{n+1}\|$, $n=0,1,\ldots$. Here, $\mu^n=\mu(\phi^n)=E(\phi^n)+\frac{\beta}{2}\int_{\mathbb{R}^d} |\phi^n|^4\mathrm{d}\mathbf{x}$ and $\vartheta^n\geq0$ serves as a stabilization factor which can be properly designed to such that the time step $\tau>0$ can be chosen as large as possible \cite{BC2013KRM,BD2004SISC,LC2021SISC}.

\footnotesize
\bibliographystyle{abbrv}
\bibliography{references}

\end{document}